%% file: series.tex
\documentclass[12pt]{amsart}
\usepackage{
	amssymb,
	hyperref,
	mathtools,
	array,
	booktabs,
	longtable
}

\allowdisplaybreaks

\providecommand{\nv}{\mathbf{n}}
\providecommand{\casimir}{\mathbf{c}}
\providecommand{\DataUnit}{\mathbf{1}}
\providecommand{\DataV}[4]{V_{#1}(#2,#3,#4)}

\providecommand{\SL}{\operatorname{SL}}
\providecommand{\GL}{\operatorname{GL}}
\providecommand{\oO}{\operatorname{O}}
\providecommand{\Aut}{\operatorname{Aut}}
\providecommand{\fS}{\mathfrak{S}}
\newcommand{\largewedge}{\mbox{\Large $\wedge$}}
\newcommand\bN{\mathbb{N}}
\newcommand\bQ{\mathbb{Q}}

\begin{document}
\title{The exceptional series through the fifth tensor power}
\author{Bruce Westbury}
\date{\today}
\begin{abstract}
We give a uniform structure for the first five tensor powers of the adjoint representation
of an exceptional simple Lie algebra.
\end{abstract}
\maketitle
\tableofcontents

\section{Introduction}
The proposal that the exceptional simple Lie algebras should form a series originated in a 1995 preprint by Pierre Vogel which was eventually published in \cite{Vogel2011}. The motivation was to study Vassiliev invariants. The popular take away is that this paper assigns a point in the projective plane to each simple Lie algebra. Let $L$ be a simple Lie algebra (other than $\SL(2)$). Write
the exterior square of $L$ as $L\oplus X$ and the symmetric square as $I\oplus Y$. The key result
is that $X$ is irreducible and $Y$ is the sum of two or three irreducible representations. Vogel then observed that the simple Lie algebras with $Y$ the sum of two irreducible representations lie on a 
line in the Vogel plane and that all five exceptional simple Lie algebras have this property and so
give points on this line. This line is now known as the exceptional series. Our choice of
coordinate is $\nv$ and we take the exceptional series to be the following eight points
with associated coordinate
\begin{longtable}{ccccccccc}
	\toprule
	& A1 & A2 & G2 & D4 & F4 & E6 & E7 & E8 \\ 
	\midrule\endhead
	$\nv$ & 1/3  & 1/2 & 2/3 & 1 & 3/2 & 2 & 3 & 5 \\
	\bottomrule
\end{longtable}
We omit the trivial group at $\nv=1/5$ as it is uninteresting. 
We have not included the supergroup at $\nv=1/4$ and the intermediate Lie algebra at $\nv=4$.

The discovery of the exceptional series pointed to the remarkable idea that the tensor powers of
the adjoint representation of the Lie algebras in the exceptional series exhibit common structure. 
This idea was developed in \cite{deligne1996b}, \cite{deligne1996a}.

This paper is a successor to \cite{cohen1996} which tabulates the structure common to the $k$-th 
tensor powers of the adjoint representations for $0\leqslant k \leqslant 4$. The main purpose
of this paper to extend these tables to the fifth tensor power. We also make the following
improvements:
\begin{itemize}
	\item We introduce systematic labels instead of the ad hoc letters in \cite{cohen1996}.
	\item We tabulate the quantum dimensions instead of the classical dimensions.
	\item We regard the orthogonal functions as the basic polynomial functors instead
	of the Schur functions.
\end{itemize}

The tables are split into two sections. The first section is the generic structure and consists of
Casimir values, branching rules, quantum dimensions and polynomial functors. There is an unexpected
involution on the generic structure. This is given by an involution on the labels together
with corresponding involutions on the Casimir values and quantum dimensions.
The second section gives the character data for each entry.

This paper is supported by \cite{Westbury2026}. The purpose of this repository is to make the
data available to any computer. The data is written into json files so there is no dependency
on operating system or software. The depository also contains two audits of the data, one
using Sage and one written in C. The Sage audit omits an audit of the polynomial functors.
The C audit omits any audit that involves arithmetic with rational functions. There is the
option to include these audits if \href{https://flintlib.org/}{flint} is available.

The audits can be described succinctly in terms of $\lambda$-rings. The branching rules
define a commutative graded ring (up to degree five) with a distinguished positive cone.
The $\lambda$-ring structure is given by the Schur functors and the orthogonal functors
preserve the positive cone. The ring of rational functions $\bQ(Q,q)$ has Adams operations
given by $\psi_k\colon Q\mapsto Q^k, q\mapsto q^k$ and these give $\bQ(Q,q)$ the structure
of a $\lambda$-ring. The quantum dimension is a morphism of $\lambda$-rings. The involution
on labels is a ring involution which preserves the positive cone (but is not an involution
of $\lambda$-rings) and there is an involution on $\bQ(Q,q)$ compatible with the quantum dimension.
For each point, $L$, in the exceptional series we have 
a morphism of $\lambda$-rings to the representation ring of $\Aut(L)$ (which does not preserve
the positive cone). The generic quantum dimensions interpolate the quantum dimensions
in the following sense. The $q$-deformation of $\Aut(L)$ has the same representation ring
as $\Aut(L)$ but the quantum dimension is a morphism of $\lambda$-rings to $\bQ(q)$.
Let $n=r/s$ be the coordinate of $L$. There is a homomorphism of $\lambda$-rings between
Laurent polynomial rings $\bQ[Q^{\pm 1},q^{\pm 1}] \to \bQ[q^{\pm 1}]$ given by $Q\mapsto Q^r$, $q\mapsto q^s$. This homomorphism does not extend to $\bQ(Q,q)$ but is defined on the generic
quantum dimensions and maps each generic quantum dimension to the corresponding quantum dimension.
Similarly the generic Casimirs interpolate the Casimir values of the corresponding representations.

\subsection{Prior information}
There are three results that were used in the discovery of the data of this paper
which are not audited.
\subsubsection{Cubic Hecke algebras}\label{hecke}
For $n>1$, let $B_n^{(3)}$ be the quotient of the braid group
$B_n$ by the relations $\sigma_i^3=1$ for $1\leqslant i\leqslant n-1$.
It is known, see \cite{coxeter1957} and \cite{assion1978},
that $B_n^{(3)}$ is finite for $n=2,3,4,5$ and is infinite for $n\geqslant 6$. The groups
$B_n^{(3)}$ for $n=3,4,5$ are the complex reflection groups with Shephard-Todd numbers 4, 25, 32.
respectively. The group $B_3^{(3)}$ has order 24 and the dimensions of the irreducible representations are 3,2,2,2,1,1,1.

We choose one representation of dimension 2, $U$. The matrices $H$ are the branching
rules for the representations of $B_3^{(k)}$, $0\leqslant k\leqslant 2$ and the
representations of $B_3^{(k)}$, $3\leqslant k\leqslant 5$ whose restriction to $B_3^{(3)}$
does not contain $U$. This is independent of the choice of $U$. It follows from \cite{marin2012}
that these are the branching rules for the algebras $Q_n$ in \cite{marin2022}.

The fact that $B_6^{(3)}$ is infinite blocks us from continuing our investigations to the
sixth tensor powers of the adjoint representations.
 
\subsubsection{Laurent-Manivel cone}
The Laurent-Manivel cone was introduced in \cite{Landsberg2002}. This is a cone in the dominant
weights of each the simple 

This gives the quantum dimension and Casimir for all labels in the Landsberg-Manivel cone and their partners.

\subsubsection{Universal characters}
A sequence of universal characters was studied in \cite{Westbury2006a}. For each Lie algebra $L$, there are characters $X_k(L)$ for $k\geqslant 1$. For small $k$, this is characterised as the 
eigenspace of the action of the Casimir, $\casimir$, on the $k$-th exterior power of $L$. The eigenvalue is $6k\nv$ with the normalisation that $\casimir(L) = 6 \nv$. The label for $X_k$ in the tables is $V_k(0,1,0)$. Each $X_k(L)$ is given by applying a polynomial functor to $L$.
For $k\leqslant 5$, the polynomial functors are, from \cite[(8)]{Westbury2006a},

\begin{align*}
X_1(L) &= L \\
X_2(L) &= \largewedge^2(L) - L \\
X_3(L) &= \largewedge^3(L) - L\otimes L + L \\
X_4(L) &= \largewedge^4(L) - \largewedge^2(L)\otimes L + \largewedge^2(L) + L\otimes L - L \\
X_5(L) &= \largewedge^5(L) - \largewedge^3(L)\otimes L + 2\, \largewedge^2(L)\otimes L - 2\, L\otimes L + L \\
\end{align*}

This gives the quantum dimension of $V_k(0,1,0)$ for all $k$, and the character of $V_k(0,1,0)$ for all members of the exceptional series and all $k$.

\section{Generic structure}
The generic labels at level $k$ for $k=2,3,4,5$ are the columns of the matrix 
$H_2H_3\dotsb H_k$. If $[a,b,c]$ is a column at level $k$ then we denote the 
representation by $V_k(a,b,c)$. This implies that the dimension of $V_k(a,b,c)$ is $a+b+c$.
The involution on labels is $V_k(a,b,c)\mapsto V_k(c,b,a)$.

\subsection{Labels and Casimirs}
The second column gives the name in \cite{cohen1996} and is included for convenience.
The third column gives the Casimir. This is normalised so that the Casimir of $L$ is $6\,\nv$.
The fourth column gives the Landsberg-Manivel coordinates if the label is in the
Landsberg-Manivel cone and is blank otherwise. This is also included for convenience.

The involution on the Casimir is $r\,\nv+s \mapsto (r-s)\nv-s$. This means that if
$V_k(a,b,c)$ has Casimir $r\,\nv+s$ then $V_k(c,b,a)$ has Casimir $(r-s)\nv-s$.

\input{labels/level0.tex}
\input{labels/level1.tex}
\input{labels/level2.tex}
\input{labels/level3.tex}
\input{labels/level4.tex}

\input{labels/level5.tex}
\subsection{Bratteli data}
The matrices representing the branching rules are built from matrices $H_k$ and $A_k$.
The method for calculating the matrices $H_k$ is given in \S\ref{hecke}. In this project
the matrices $H_k$ and $A_k$ are taken as initial data and are not audited.

The branching matrices are given by a tridiagonal block matrix construction with the
matrices $H_k$ on the super-diagonal, the matrices $A_k$ on the diagonal
and the transposes of the matrices $H_k$ on the sub-diagonal.

We given the branching rules from level 4 to level 5. The lower branching rules
are given by dropping final rows and final columns.
The column sizes are 1,1,3,6,15,30 and the row sizes are 1,3,6,15,30.
The matrix entries are:

\begin{equation*}
\begin{pmatrix}
0 & H_0 & 0 & 0 & 0 & 0 \\
H_0^t & A_1 & H_1 & 0 & 0 & 0 \\
0 & H_1^t & A_2 & H_2 & 0 & 0 \\
0 & 0 & H_2^t & A_3 & H_3 & 0 \\
0 & 0 & 0 & H_3^t & A_4 & H_4
\end{pmatrix}
\end{equation*}
These matrices are present in the archive but are omitted here.

\input{bratteli/H_0.tex}
\input{bratteli/H_1.tex}
\input{bratteli/H_2.tex}
\input{bratteli/H_3.tex}
\input{bratteli/H_4.tex}
\input{bratteli/A_2.tex}
\input{bratteli/A_3.tex}
\input{bratteli/A_4.tex}

\subsection{Product decompositions}
This table gives the product of a level 2 label and a level 2 label.
\input{products/level2_level2_products.tex}
This table gives the product of a level 2 label and a level 3 label.
\input{products/level2_level3_products.tex}

\subsection{Quantum dimensions}
Each dimension formula is given by two lists of quantum integers, $N(V)$ and $D(V)$,
of the same length, see \cite[\S1.2]{morrison2024}. The dimension is $N(V)/D(V)$ but the advantage of this notation is that there are several interpretations.

The two-variable quantum integer interpretation is 
\[ [r\nv+s] = Q^rq^s - Q^{-r}q^{-s} \]
This interpretation gives the dimension as an element of $\bQ(Q,q)$.
The classical version is the interpretation
\[ [r\nv+s] = (r\nv+s) \]
which gives the dimension as an element of $\bQ(\nv)$. For level at most 4 these are
the dimension formulae in \cite{cohen1996}.
For the point $\nv = n$ in the exceptional series,  we take the interpretation
\[ [r\nv+s] = [r\,n+s] \] where $[r\,n+s]$ is a quantum integer. This gives the dimension
as an element of $\bQ(q)$. Finally, we can take the interpretation
\[ [r\nv+s] = (r\,n+s) \]  
This gives the dimension as an element of $\bQ$.

The involution is $[r\nv+s] \mapsto [(s-r)\nv+s]$ which corresponds to $Q\mapsto Q^{-1}$, $q\mapsto Qq$.

\input{qdim/level0.tex}
\input{qdim/level1.tex}
\input{qdim/level2.tex}
\input{qdim/level3.tex}
\input{qdim/level4.tex}
\input{qdim/level5.tex}

\section{Orthogonal functors}
The orthogonal functors are the polynomial functors associated to the orthogonal
functions in \cite{Koike1987}. Let $L$ be a simple Lie algebra, then the Killing form
gives a group homomorphism $\Aut(L) \to \oO(L)$. The orthogonal functors give the 
restriction of an irreducible representation $\oO(L)$ just as the Schur functors 
give the restriction of an irreducible representation of $\GL(L)$.

For the convenience of the reader we record the change of basis between the Schur functions and the orthogonal functions.
\paragraph{Level 2}
\begin{align*}
	o_{(2)} &= s_{(2)} - s_{()} & s_{(2)} &= o_{(2)} + o_{()} \\
	o_{(1,1)} &= s_{(1,1)} & s_{(1,1)} &= o_{(1,1)} 
\end{align*}
\paragraph{Level 3}
\begin{align*}
	o_{(3)} &= s_{(3)} - s_{(1)} & s_{(3)} &= o_{(3)} + o_{(1)} \\
	o_{(2,1)} &= s_{(2,1)} - s_{(1)} & s_{(2,1)} &= o_{(2,1)} + o_{(1)} \\
	o_{(1,1,1)} &= s_{(1,1,1)} & s_{(1,1,1)} &= o_{(1,1,1)}
\end{align*}
\paragraph{Level 4}
\begin{align*}
	o_{(4)} &= s_{(4)} - s_{(2)}  & s_{(4)} &= o_{(4)}  + o_{(2)}  + o_{()} \\
	o_{(3,1)} &= s_{(3,1)} - s_{(2)} - s_{(1,1)} + s_{(1)} 
	& s_{(3,1)} &= o_{(3,1)}  + o_{(2)} + o_{(1,1)} \\
	o_{(2,2)} &= s_{(2,2)} - s_{(2)} & s_{(2,2)} &= o_{(2,2)}  + o_{(2)}  + o_{()} \\
	o_{(2,1,1)} &= s_{(2,1,1)} - s_{(1,1)}  & s_{(2,1,1)} &= o_{(2,1,1)} + o_{(1,1)} \\
	o_{(1,1,1,1)} &= s_{(1,1,1,1)} & s_{(1,1,1,1)} &= o_{(1,1,1,1)}
\end{align*}
\paragraph{Level 5}
\begin{align*}
	o_{(5)} &= s_{(5)} - s_{(3)} 
	& s_{(5)} &= o_{(5)} + o_{(3)}  + o_{(1)} \\
	o_{(4,1)} &= s_{(4,1)} - s_{(3)} - s_{(2,1)} + s_{(1)} 
	& s_{(4,1)} &= o_{(4,1)} + o_{(3)} + o_{(2,1)}  + o_{(1)}  \\
	o_{(3,2)} &= s_{(3,2)} - s_{(3)} - s_{(2,1)} + s_{(1)} 
	& s_{(3,2)} &= o_{(3,2)} + o_{(3)} + o_{(2,1)}  + o_{(1)} \\
	o_{(3,1,1)} &= s_{(3,1,1)} - s_{(2,1)} - s_{(1,1,1)} + s_{(1)}
	& s_{(3,1,1)} &= o_{(3,1,1)} + o_{(2,1)}  + o_{(1,1,1)} \\
	o_{(2,2,1)} &= s_{(2,2,1)} - s_{(2,1)} 
	& s_{(2,2,1)} &= o_{(2,2,1)} + o_{(2,1)}  + o_{(1)} \\
	o_{(2,1,1,1)} &= s_{(2,1,1,1)} - s_{(1,1,1)} 
	& s_{(2,1,1,1)} &= o_{(2,1,1,1)} + o_{(1,1,1)} \\
	o_{(1,1,1,1,1)} &= s_{(1,1,1,1,1)} 
	& s_{(1,1,1,1,1)} &= o_{(1,1,1,1,1)}
\end{align*}

\input{orthogonal/degree1.tex}
\input{orthogonal/degree2.tex}
\input{orthogonal/degree3.tex}
\input{orthogonal/degree4.tex}
\input{orthogonal/degree5.tex}

\subsection{Symmetric and exterior squares}
This section gives the symmetric and exterior squares of the level 2 labels.
\input{products/level2_square_functors.tex}

\section{Characters}
For each Cartan type in the exceptional series let $L$ be the associated simple Lie algebra.  Define $G\coloneqq \Aut(L)$, the automorphism group of the adjoint representation.
The identity component, $G_0$, is the connected algebraic group of adjoint type whose Lie algebra is $L$ and $G\cong G_0 \rtimes \Gamma$ where $\Gamma$ is the group of diagram automorphisms.

By Clifford theory (the representation theory of semidirect products), the irreducible
representations of $\Aut(L)$ are indexed by pairs consisting of an orbit of the action
of $\Gamma$ on dominant weights and an irreducible representation of the stabiliser subgroup.
In the tables we give only give the character of $G_0$. The characters of the stabilisers
up to level 4 are given in \cite{cohen1996}.

For $A2$ and $E6$, $\Gamma$ has order 2. An irreducible representation are indexed by
either an orbit of size 2 or an orbit of size 1 and a sign.
For $D4$, $\Gamma$ is the permutation group $\fS_3$.

\subsection{Characters}
\input{characters/level0.tex}
\input{characters/level1.tex}
\input{characters/level2.tex}
\input{characters/level3.tex}
\input{characters/level4.tex}
\input{characters/level5.tex}

\section{Conclusion}
The data we have presented is overdetermined; informally, there are more constraints than there are degrees of freedom. Moreover the tables and the structure interact in a complicated way. We did not
fill out each table one at a time. Instead we started with each table partially completed and incrementally filled in the tables by iteratively using one piece of the structure to determine a few
more entries in a table. We chose not to document this process as it depends on a series of choices;
the most significant being the choice of which data to assume at the start. There is also structure
that is implicit; for example, the branching rules gives linear equations but the characters are all
simple linear combinations of irreducible characters of $\Aut(L)$ with integer coefficients,
and most of the quantum dimensions can be written using the pair notation.

The fact that this data exists with all the structure is unexpected and warrants an explanation.
The open question is whether this data can be continued to the sixth, and higher, tensor powers.
A formal way to ask the question is whether the data is the first five levels of a 
commutative, graded, ring with Adams operations. The natural construction is as the Grothendieck
group of a symmetric monoidal category whose monoid of objects is $\bN$. The line pursued in
\cite{Cohen1999} and, more recently, in \cite{morrison2024},\cite{westbury2024} is that this category interpolates
the categories of invariant tensors of $L$ as a representation of $\Aut(L)$. This is not 
consistent with our data; the character $V_5(0,1,0)$ of $F4$ and $E6$ is the sum of two irreducible characters, as is the character $V_5(2,30)$ for $D4$.
The character formulae in \cite{Westbury2006a} give further examples at higher levels.

\section{Use of AI}

The software, audit infrastructure, and supporting documentation accompanying this dataset were developed with the assistance of an AI coding agent. The mathematical objectives, computational strategy, data formats, verification methodology, and interpretation of the results were conceived and directed by the author. The AI agent implemented the requested software, generated documentation, and carried out programming tasks under detailed mathematical guidance and iterative human review.

The author accepts full responsibility for the mathematical content of the project, the correctness of the published data, and the conclusions drawn from them. Confidence in the dataset does not rest on the AI-generated implementation but on independent mathematical verification. The published archive therefore contains two distinct audit implementations. The primary audit suite, written in C, performs independent calculations using LiE together with optional exact rational-function arithmetic provided by FLINT, and gives complete coverage of the archived mathematical data. A second audit suite, implemented independently in Sage, verifies all remaining mathematical structures directly from the published JSON files without using either the C implementation or LiE. The only deliberate exception is the verification of polynomial-functor calculations, since Sage does not currently provide suitable functionality for these computations. Together these independent audits provide comprehensive verification of the published dataset.

\bibliographystyle{alphaurl}
\bibliography{series}

\end{document}

%% file: labels/level0.tex
\begin{longtable}{@{} >{$}l<{$} >{$}c<{$} >{$}l<{$} >{$}c<{$} @{}}
\toprule
\text{label} & \text{Cohen--de~Man} & \casimir & \text{LM} \\
\midrule\endhead
\DataUnit & I & 0 & (0,0,0,0) \\
\bottomrule
\end{longtable}

%% file: labels/level1.tex
\begin{longtable}{@{} >{$}l<{$} >{$}c<{$} >{$}l<{$} >{$}c<{$} @{}}
\toprule
\text{label} & \text{Cohen--de~Man} & \casimir & \text{LM} \\
\midrule\endhead
L & L & 6\nv & (1,0,0,0) \\
\bottomrule
\end{longtable}

%% file: labels/level2.tex
\begin{longtable}{@{} >{$}l<{$} >{$}c<{$} >{$}l<{$} >{$}c<{$} @{}}
\toprule
\text{label} & \text{Cohen--de~Man} & \casimir & \text{LM} \\
\midrule\endhead
\DataV{2}{1}{0}{0} & Y_{2}^* & 10\nv - 2 & (0,0,0,1) \\
\DataV{2}{0}{1}{0} & X_{2} & 12\nv & (0,1,0,0) \\
\DataV{2}{0}{0}{1} & Y_{2} & 12\nv + 2 & (2,0,0,0) \\
\bottomrule
\end{longtable}

%% file: labels/level3.tex
\begin{longtable}{@{} >{$}l<{$} >{$}c<{$} >{$}l<{$} >{$}c<{$} @{}}
\toprule
\text{label} & \text{Cohen--de~Man} & \casimir & \text{LM} \\
\midrule\endhead
\DataV{3}{1}{0}{0} & Y_{3}^* & 12\nv - 6 \\
\DataV{3}{0}{1}{0} & X_{3} & 18\nv & (0,0,1,0) \\
\DataV{3}{0}{0}{1} & Y_{3} & 18\nv + 6 & (3,0,0,0) \\
\DataV{3}{0}{1}{1} & C & 18\nv + 3 & (1,1,0,0) \\
\DataV{3}{1}{1}{0} & C^* & 15\nv - 3 \\
\DataV{3}{1}{1}{1} & A & 16\nv & (1,0,0,1) \\
\bottomrule
\end{longtable}

%% file: labels/level4.tex
\begin{longtable}{@{} >{$}l<{$} >{$}c<{$} >{$}l<{$} >{$}c<{$} @{}}
\toprule
\text{label} & \text{Cohen--de~Man} & \casimir & \text{LM} \\
\midrule\endhead
\DataV{4}{1}{0}{0} & Y_{4}^* & 12\nv - 12 \\
\DataV{4}{0}{0}{1} & Y_{4} & 24\nv + 12 & (4,0,0,0,) \\
\DataV{4}{0}{1}{0} & X_{4} & 24\nv \\
\DataV{4}{0}{1}{1} & H & 24\nv + 6 & (0,2,0,0,) \\
\DataV{4}{1}{1}{0} & H^* & 18\nv - 6 \\
\DataV{4}{1}{1}{1} & J & 20\nv & (0,0,0,2,) \\
\DataV{4}{1}{2}{0} & I^* & 20\nv - 4 \\
\DataV{4}{2}{1}{0} & G^* & 16\nv - 8 \\
\DataV{4}{0}{2}{1} & I & 24\nv + 4 & (1,0,1,0,) \\
\DataV{4}{0}{1}{2} & G & 24\nv + 8 & (2,1,0,0,) \\
\DataV{4}{1}{3}{2} & F & 22\nv + 2 & (0,1,0,1,) \\
\DataV{4}{1}{2}{3} & D & 22\nv + 4 & (2,0,0,1,) \\
\DataV{4}{3}{2}{1} & D^* & 18\nv - 4 \\
\DataV{4}{2}{3}{1} & F^* & 20\nv - 2 \\
\DataV{4}{2}{4}{2} & E & 21\nv \\
\bottomrule
\end{longtable}

%% file: labels/level5.tex
\begin{longtable}{@{} >{$}l<{$} >{$}l<{$} >{$}c<{$} @{}}
\toprule
\text{label} & \casimir & \text{LM} \\
\midrule\endhead
\DataV{5}{1}{0}{0} & 10\nv - 20 \\
\DataV{5}{0}{0}{1} & 30\nv + 20 & (5,0,0,0) \\
\DataV{5}{0}{1}{0} & 30\nv \\
\DataV{5}{0}{1}{3} & 30\nv + 15 & (3,1,0,0) \\
\DataV{5}{1}{3}{0} & 25\nv - 5 \\
\DataV{5}{0}{3}{1} & 30\nv + 5 \\
\DataV{5}{3}{1}{0} & 15\nv - 15 \\
\DataV{5}{0}{2}{3} & 30\nv + 12 & (1,2,0,0) \\
\DataV{5}{0}{3}{2} & 30\nv + 8 \\
\DataV{5}{2}{3}{0} & 22\nv - 8 & (0,1,1,0) \\
\DataV{5}{3}{2}{0} & 18\nv - 12 \\
\DataV{5}{0}{3}{3} & 30\nv + 10 \\
\DataV{5}{3}{3}{0} & 20\nv - 10 & (2,0,1,0) \\
\DataV{5}{1}{6}{3} & 28\nv + 4 \\
\DataV{5}{3}{6}{1} & 24\nv - 4 & (0,0,1,1) \\
\DataV{5}{1}{3}{6} & 28\nv + 10 \\
\DataV{5}{6}{3}{1} & 18\nv - 10 & (3,0,0,1) \\
\DataV{5}{3}{6}{6} & 26\nv + 4 \\
\DataV{5}{6}{6}{3} & 22\nv - 4 & (1,0,0,2) \\
\DataV{5}{3}{9}{3} & 26\nv \\
\DataV{5}{2}{9}{9} & 28\nv + 7 \\
\DataV{5}{6}{8}{6} & 24\nv \\
\DataV{5}{3}{9}{8} & 27\nv + 5 \\
\DataV{5}{8}{9}{3} & 22\nv - 5 \\
\DataV{5}{3}{11}{6} & 27\nv + 3 \\
\DataV{5}{9}{9}{2} & 21\nv - 7 & (1,1,0,1) \\
\DataV{5}{6}{11}{3} & 24\nv - 3 \\
\DataV{5}{6}{12}{6} & 25\nv \\
\DataV{5}{6}{15}{9} & 26\nv + 2 \\
\DataV{5}{9}{15}{6} & 24\nv - 2 \\
\bottomrule
\end{longtable}

%% file: bratteli/H_0.tex
The matrix $H_0$ has shape $1\times1$.  Its nonzero entries are:
\begin{align*}
H_0\!\left(\DataUnit\right) &= \DataUnit \\
\end{align*}

%% file: bratteli/H_1.tex
The matrix $H_1$ has shape $1\times3$.  Its nonzero entries are:
\begin{align*}
H_1\!\left(L\right) &= \DataV{2}{1}{0}{0} + \DataV{2}{0}{1}{0} + \DataV{2}{0}{0}{1} \\
\end{align*}

%% file: bratteli/H_2.tex
The matrix $H_2$ has shape $3\times6$.  Its nonzero entries are:
\begin{align*}
H_2\!\left(\DataV{2}{1}{0}{0}\right) &= \DataV{3}{1}{0}{0} + \DataV{3}{1}{1}{0} + \DataV{3}{1}{1}{1} \\
H_2\!\left(\DataV{2}{0}{1}{0}\right) &= \DataV{3}{0}{1}{0} + \DataV{3}{0}{1}{1} + \DataV{3}{1}{1}{0} + \DataV{3}{1}{1}{1} \\
H_2\!\left(\DataV{2}{0}{0}{1}\right) &= \DataV{3}{0}{0}{1} + \DataV{3}{0}{1}{1} + \DataV{3}{1}{1}{1} \\
\end{align*}

%% file: bratteli/H_3.tex
The matrix $H_3$ has shape $6\times15$.  Its nonzero entries are:
\begin{align*}
H_3\!\left(\DataV{3}{1}{0}{0}\right) &= \DataV{4}{1}{0}{0} + \DataV{4}{2}{1}{0} + \DataV{4}{3}{2}{1} \\
H_3\!\left(\DataV{3}{0}{1}{0}\right) &= \DataV{4}{0}{1}{0} + \DataV{4}{1}{2}{0} + \DataV{4}{0}{2}{1} + \DataV{4}{1}{3}{2} \\
&\quad {}+ \DataV{4}{2}{3}{1} + \DataV{4}{2}{4}{2} \\
H_3\!\left(\DataV{3}{0}{0}{1}\right) &= \DataV{4}{0}{0}{1} + \DataV{4}{0}{1}{2} + \DataV{4}{1}{2}{3} \\
H_3\!\left(\DataV{3}{0}{1}{1}\right) &= \DataV{4}{0}{1}{1} + \DataV{4}{0}{2}{1} + \DataV{4}{0}{1}{2} + \DataV{4}{1}{3}{2} \\
&\quad {}+ \DataV{4}{1}{2}{3} + \DataV{4}{2}{4}{2} \\
H_3\!\left(\DataV{3}{1}{1}{0}\right) &= \DataV{4}{1}{1}{0} + \DataV{4}{1}{2}{0} + \DataV{4}{2}{1}{0} + \DataV{4}{3}{2}{1} \\
&\quad {}+ \DataV{4}{2}{3}{1} + \DataV{4}{2}{4}{2} \\
H_3\!\left(\DataV{3}{1}{1}{1}\right) &= \DataV{4}{1}{1}{1} + \DataV{4}{1}{3}{2} + \DataV{4}{1}{2}{3} + \DataV{4}{3}{2}{1} \\
&\quad {}+ \DataV{4}{2}{3}{1} + \DataV{4}{2}{4}{2} \\
\end{align*}

%% file: bratteli/H_4.tex
The matrix $H_4$ has shape $15\times30$.  Its nonzero entries are:
\begin{align*}
H_4\!\left(\DataV{4}{1}{0}{0}\right) &= \DataV{5}{1}{0}{0} + \DataV{5}{3}{1}{0} + \DataV{5}{6}{3}{1} \\
H_4\!\left(\DataV{4}{0}{0}{1}\right) &= \DataV{5}{0}{0}{1} + \DataV{5}{0}{1}{3} + \DataV{5}{1}{3}{6} \\
H_4\!\left(\DataV{4}{0}{1}{0}\right) &= \DataV{5}{0}{1}{0} + \DataV{5}{1}{3}{0} + \DataV{5}{0}{3}{1} + \DataV{5}{1}{6}{3} \\
&\quad {}+ \DataV{5}{3}{6}{1} + \DataV{5}{3}{9}{3} + \DataV{5}{3}{11}{6} + \DataV{5}{6}{11}{3} \\
&\quad {}+ \DataV{5}{6}{12}{6} + \DataV{5}{6}{15}{9} + \DataV{5}{9}{15}{6} \\
H_4\!\left(\DataV{4}{0}{1}{1}\right) &= \DataV{5}{0}{2}{3} + \DataV{5}{0}{3}{2} + \DataV{5}{2}{9}{9} + \DataV{5}{3}{11}{6} \\
H_4\!\left(\DataV{4}{1}{1}{0}\right) &= \DataV{5}{2}{3}{0} + \DataV{5}{3}{2}{0} + \DataV{5}{9}{9}{2} + \DataV{5}{6}{11}{3} \\
H_4\!\left(\DataV{4}{1}{1}{1}\right) &= \DataV{5}{3}{6}{6} + \DataV{5}{6}{6}{3} + \DataV{5}{6}{12}{6} \\
H_4\!\left(\DataV{4}{1}{2}{0}\right) &= \DataV{5}{1}{3}{0} + \DataV{5}{2}{3}{0} + \DataV{5}{3}{3}{0} + \DataV{5}{3}{6}{1} \\
&\quad {}+ \DataV{5}{3}{9}{3} + \DataV{5}{8}{9}{3} + \DataV{5}{9}{9}{2} + \DataV{5}{6}{11}{3} \\
&\quad {}+ \DataV{5}{9}{15}{6} \\
H_4\!\left(\DataV{4}{2}{1}{0}\right) &= \DataV{5}{3}{1}{0} + \DataV{5}{3}{2}{0} + \DataV{5}{3}{3}{0} + \DataV{5}{6}{3}{1} \\
&\quad {}+ \DataV{5}{8}{9}{3} + \DataV{5}{9}{9}{2} \\
H_4\!\left(\DataV{4}{0}{2}{1}\right) &= \DataV{5}{0}{3}{1} + \DataV{5}{0}{3}{2} + \DataV{5}{0}{3}{3} + \DataV{5}{1}{6}{3} \\
&\quad {}+ \DataV{5}{3}{9}{3} + \DataV{5}{2}{9}{9} + \DataV{5}{3}{9}{8} + \DataV{5}{3}{11}{6} \\
&\quad {}+ \DataV{5}{6}{15}{9} \\
H_4\!\left(\DataV{4}{0}{1}{2}\right) &= \DataV{5}{0}{1}{3} + \DataV{5}{0}{2}{3} + \DataV{5}{0}{3}{3} + \DataV{5}{1}{3}{6} \\
&\quad {}+ \DataV{5}{2}{9}{9} + \DataV{5}{3}{9}{8} \\
H_4\!\left(\DataV{4}{1}{3}{2}\right) &= \DataV{5}{1}{6}{3} + \DataV{5}{3}{6}{6} + \DataV{5}{2}{9}{9} + \DataV{5}{3}{11}{6} \\
&\quad {}+ \DataV{5}{6}{12}{6} + \DataV{5}{6}{15}{9} + \DataV{5}{9}{15}{6} \\
H_4\!\left(\DataV{4}{1}{2}{3}\right) &= \DataV{5}{1}{3}{6} + \DataV{5}{3}{6}{6} + \DataV{5}{2}{9}{9} + \DataV{5}{6}{8}{6} \\
&\quad {}+ \DataV{5}{3}{9}{8} + \DataV{5}{6}{15}{9} \\
H_4\!\left(\DataV{4}{3}{2}{1}\right) &= \DataV{5}{6}{3}{1} + \DataV{5}{6}{6}{3} + \DataV{5}{6}{8}{6} + \DataV{5}{8}{9}{3} \\
&\quad {}+ \DataV{5}{9}{9}{2} + \DataV{5}{9}{15}{6} \\
H_4\!\left(\DataV{4}{2}{3}{1}\right) &= \DataV{5}{3}{6}{1} + \DataV{5}{6}{6}{3} + \DataV{5}{9}{9}{2} + \DataV{5}{6}{11}{3} \\
&\quad {}+ \DataV{5}{6}{12}{6} + \DataV{5}{6}{15}{9} + \DataV{5}{9}{15}{6} \\
H_4\!\left(\DataV{4}{2}{4}{2}\right) &= \DataV{5}{3}{9}{3} + \DataV{5}{6}{8}{6} + \DataV{5}{3}{9}{8} + \DataV{5}{8}{9}{3} \\
&\quad {}+ \DataV{5}{3}{11}{6} + \DataV{5}{6}{11}{3} + \DataV{5}{6}{12}{6} + \DataV{5}{6}{15}{9} \\
&\quad {}+ \DataV{5}{9}{15}{6} \\
\end{align*}

%% file: bratteli/A_2.tex
The matrix $A_2$ has shape $3\times3$.  Its nonzero entries are:
\begin{align*}
A_2\!\left(\DataV{2}{1}{0}{0}\right) &= \DataV{2}{1}{0}{0} + \DataV{2}{0}{1}{0} \\
A_2\!\left(\DataV{2}{0}{1}{0}\right) &= \DataV{2}{1}{0}{0} + \DataV{2}{0}{1}{0} + \DataV{2}{0}{0}{1} \\
A_2\!\left(\DataV{2}{0}{0}{1}\right) &= \DataV{2}{0}{1}{0} + \DataV{2}{0}{0}{1} \\
\end{align*}

%% file: bratteli/A_3.tex
The matrix $A_3$ has shape $6\times6$.  Its nonzero entries are:
\begin{align*}
A_3\!\left(\DataV{3}{1}{0}{0}\right) &= \DataV{3}{1}{0}{0} + \DataV{3}{1}{1}{0} \\
A_3\!\left(\DataV{3}{0}{1}{0}\right) &= \DataV{3}{0}{1}{0} + \DataV{3}{0}{1}{1} + \DataV{3}{1}{1}{0} + \DataV{3}{1}{1}{1} \\
A_3\!\left(\DataV{3}{0}{0}{1}\right) &= \DataV{3}{0}{0}{1} + \DataV{3}{0}{1}{1} \\
A_3\!\left(\DataV{3}{0}{1}{1}\right) &= \DataV{3}{0}{1}{0} + \DataV{3}{0}{0}{1} + 2\DataV{3}{0}{1}{1} + \DataV{3}{1}{1}{1} \\
A_3\!\left(\DataV{3}{1}{1}{0}\right) &= \DataV{3}{1}{0}{0} + \DataV{3}{0}{1}{0} + 2\DataV{3}{1}{1}{0} + \DataV{3}{1}{1}{1} \\
A_3\!\left(\DataV{3}{1}{1}{1}\right) &= \DataV{3}{0}{1}{0} + \DataV{3}{0}{1}{1} + \DataV{3}{1}{1}{0} + 2\DataV{3}{1}{1}{1} \\
\end{align*}

%% file: bratteli/A_4.tex
The matrix $A_4$ has shape $15\times15$.  Its nonzero entries are:
\begin{align*}
A_4\!\left(\DataV{4}{1}{0}{0}\right) &= \DataV{4}{1}{0}{0} + \DataV{4}{2}{1}{0} \\
A_4\!\left(\DataV{4}{0}{0}{1}\right) &= \DataV{4}{0}{0}{1} + \DataV{4}{0}{1}{2} \\
A_4\!\left(\DataV{4}{0}{1}{0}\right) &= 2\DataV{4}{0}{1}{0} + \DataV{4}{1}{2}{0} + \DataV{4}{0}{2}{1} + \DataV{4}{1}{3}{2} \\
&\quad {}+ \DataV{4}{2}{3}{1} + \DataV{4}{2}{4}{2} \\
A_4\!\left(\DataV{4}{0}{1}{1}\right) &= \DataV{4}{0}{1}{1} + \DataV{4}{0}{2}{1} + \DataV{4}{0}{1}{2} + \DataV{4}{1}{3}{2} \\
A_4\!\left(\DataV{4}{1}{1}{0}\right) &= \DataV{4}{1}{1}{0} + \DataV{4}{1}{2}{0} + \DataV{4}{2}{1}{0} + \DataV{4}{2}{3}{1} \\
A_4\!\left(\DataV{4}{1}{1}{1}\right) &= \DataV{4}{1}{1}{1} + \DataV{4}{1}{3}{2} + \DataV{4}{2}{3}{1} \\
A_4\!\left(\DataV{4}{1}{2}{0}\right) &= \DataV{4}{0}{1}{0} + \DataV{4}{1}{1}{0} + 2\DataV{4}{1}{2}{0} + \DataV{4}{2}{1}{0} \\
&\quad {}+ \DataV{4}{3}{2}{1} + \DataV{4}{2}{3}{1} + \DataV{4}{2}{4}{2} \\
A_4\!\left(\DataV{4}{2}{1}{0}\right) &= \DataV{4}{1}{0}{0} + \DataV{4}{1}{1}{0} + \DataV{4}{1}{2}{0} + 2\DataV{4}{2}{1}{0} \\
&\quad {}+ \DataV{4}{3}{2}{1} \\
A_4\!\left(\DataV{4}{0}{2}{1}\right) &= \DataV{4}{0}{1}{0} + \DataV{4}{0}{1}{1} + 2\DataV{4}{0}{2}{1} + \DataV{4}{0}{1}{2} \\
&\quad {}+ \DataV{4}{1}{3}{2} + \DataV{4}{1}{2}{3} + \DataV{4}{2}{4}{2} \\
A_4\!\left(\DataV{4}{0}{1}{2}\right) &= \DataV{4}{0}{0}{1} + \DataV{4}{0}{1}{1} + \DataV{4}{0}{2}{1} + 2\DataV{4}{0}{1}{2} \\
&\quad {}+ \DataV{4}{1}{2}{3} \\
A_4\!\left(\DataV{4}{1}{3}{2}\right) &= \DataV{4}{0}{1}{0} + \DataV{4}{0}{1}{1} + \DataV{4}{1}{1}{1} + \DataV{4}{0}{2}{1} \\
&\quad {}+ 2\DataV{4}{1}{3}{2} + \DataV{4}{1}{2}{3} + \DataV{4}{2}{3}{1} + \DataV{4}{2}{4}{2} \\
A_4\!\left(\DataV{4}{1}{2}{3}\right) &= \DataV{4}{0}{2}{1} + \DataV{4}{0}{1}{2} + \DataV{4}{1}{3}{2} + 2\DataV{4}{1}{2}{3} \\
&\quad {}+ \DataV{4}{2}{4}{2} \\
A_4\!\left(\DataV{4}{3}{2}{1}\right) &= \DataV{4}{1}{2}{0} + \DataV{4}{2}{1}{0} + 2\DataV{4}{3}{2}{1} + \DataV{4}{2}{3}{1} \\
&\quad {}+ \DataV{4}{2}{4}{2} \\
A_4\!\left(\DataV{4}{2}{3}{1}\right) &= \DataV{4}{0}{1}{0} + \DataV{4}{1}{1}{0} + \DataV{4}{1}{1}{1} + \DataV{4}{1}{2}{0} \\
&\quad {}+ \DataV{4}{1}{3}{2} + \DataV{4}{3}{2}{1} + 2\DataV{4}{2}{3}{1} + \DataV{4}{2}{4}{2} \\
A_4\!\left(\DataV{4}{2}{4}{2}\right) &= \DataV{4}{0}{1}{0} + \DataV{4}{1}{2}{0} + \DataV{4}{0}{2}{1} + \DataV{4}{1}{3}{2} \\
&\quad {}+ \DataV{4}{1}{2}{3} + \DataV{4}{3}{2}{1} + \DataV{4}{2}{3}{1} + 3\DataV{4}{2}{4}{2} \\
\end{align*}

%% file: products/level2_level2_products.tex
\subsubsection*{Level-two products}
\begin{align*}
\DataV{2}{1}{0}{0}\otimes\DataV{2}{1}{0}{0} &= \DataUnit + \DataV{4}{1}{0}{0} + 2\DataV{2}{1}{0}{0} + \DataV{4}{1}{1}{0} \\
&\quad {}+ \DataV{4}{1}{1}{1} + \DataV{3}{1}{0}{0} + 2\DataV{3}{1}{1}{0} + \DataV{4}{2}{3}{1} \\
&\quad {}+ \DataV{3}{0}{1}{0} + \DataV{2}{0}{1}{0} + L + \DataV{4}{2}{1}{0} \\
&\quad {}+ \DataV{3}{1}{1}{1} + \DataV{4}{3}{2}{1} + \DataV{2}{0}{0}{1} \\
\DataV{2}{1}{0}{0}\otimes\DataV{2}{0}{1}{0} &= \DataV{2}{1}{0}{0} + \DataV{3}{1}{0}{0} + 2\DataV{3}{1}{1}{0} + \DataV{4}{2}{3}{1} \\
&\quad {}+ \DataV{3}{0}{1}{0} + 2\DataV{2}{0}{1}{0} + \DataV{4}{1}{2}{0} + \DataV{4}{1}{3}{2} \\
&\quad {}+ L + \DataV{4}{2}{1}{0} + 2\DataV{3}{1}{1}{1} + \DataV{4}{3}{2}{1} \\
&\quad {}+ \DataV{4}{2}{4}{2} + \DataV{3}{0}{1}{1} + \DataV{2}{0}{0}{1} \\
\DataV{2}{1}{0}{0}\otimes\DataV{2}{0}{0}{1} &= \DataV{2}{1}{0}{0} + \DataV{3}{1}{1}{0} + \DataV{3}{0}{1}{0} + \DataV{2}{0}{1}{0} \\
&\quad {}+ \DataV{3}{1}{1}{1} + \DataV{4}{3}{2}{1} + \DataV{4}{2}{4}{2} + \DataV{3}{0}{1}{1} \\
&\quad {}+ \DataV{2}{0}{0}{1} + \DataV{4}{1}{2}{3} \\
\DataV{2}{0}{1}{0}\otimes\DataV{2}{0}{1}{0} &= \DataUnit + 2\DataV{2}{1}{0}{0} + \DataV{4}{1}{1}{0} + \DataV{4}{1}{1}{1} \\
&\quad {}+ \DataV{3}{1}{0}{0} + 2\DataV{3}{1}{1}{0} + \DataV{4}{2}{3}{1} + 2\DataV{3}{0}{1}{0} \\
&\quad {}+ \DataV{4}{0}{1}{0} + 2\DataV{2}{0}{1}{0} + \DataV{4}{1}{2}{0} + \DataV{4}{1}{3}{2} \\
&\quad {}+ \DataV{4}{0}{1}{1} + L + 3\DataV{3}{1}{1}{1} + \DataV{4}{3}{2}{1} \\
&\quad {}+ 2\DataV{4}{2}{4}{2} + \DataV{4}{0}{2}{1} + 2\DataV{3}{0}{1}{1} + 2\DataV{2}{0}{0}{1} \\
&\quad {}+ \DataV{4}{1}{2}{3} + \DataV{3}{0}{0}{1} \\
\DataV{2}{0}{1}{0}\otimes\DataV{2}{0}{0}{1} &= \DataV{2}{1}{0}{0} + \DataV{3}{1}{1}{0} + \DataV{4}{2}{3}{1} + \DataV{3}{0}{1}{0} \\
&\quad {}+ 2\DataV{2}{0}{1}{0} + \DataV{4}{1}{3}{2} + L + 2\DataV{3}{1}{1}{1} \\
&\quad {}+ \DataV{4}{2}{4}{2} + \DataV{4}{0}{2}{1} + 2\DataV{3}{0}{1}{1} + \DataV{2}{0}{0}{1} \\
&\quad {}+ \DataV{4}{1}{2}{3} + \DataV{4}{0}{1}{2} + \DataV{3}{0}{0}{1} \\
\DataV{2}{0}{0}{1}\otimes\DataV{2}{0}{0}{1} &= \DataUnit + \DataV{2}{1}{0}{0} + \DataV{4}{1}{1}{1} + \DataV{3}{0}{1}{0} \\
&\quad {}+ \DataV{2}{0}{1}{0} + \DataV{4}{1}{3}{2} + \DataV{4}{0}{1}{1} + L \\
&\quad {}+ \DataV{3}{1}{1}{1} + 2\DataV{3}{0}{1}{1} + 2\DataV{2}{0}{0}{1} + \DataV{4}{1}{2}{3} \\
&\quad {}+ \DataV{4}{0}{1}{2} + \DataV{3}{0}{0}{1} + \DataV{4}{0}{0}{1} \\
\end{align*}

%% file: products/level2_level3_products.tex
\subsubsection*{Level-two by level-three products}
\begin{align*}
\DataV{2}{1}{0}{0}\otimes\DataV{3}{1}{0}{0} &= \DataV{4}{1}{0}{0} + \DataV{2}{1}{0}{0} + \DataV{4}{1}{1}{0} + 2\DataV{3}{1}{0}{0} \\
&\quad {}+ 2\DataV{3}{1}{1}{0} + \DataV{4}{2}{3}{1} + \DataV{2}{0}{1}{0} + \DataV{4}{1}{2}{0} \\
&\quad {}+ L + 2\DataV{4}{2}{1}{0} + \DataV{3}{1}{1}{1} + \DataV{4}{3}{2}{1} \\
&\quad {}+ \DataV{5}{1}{0}{0} + \DataV{5}{3}{1}{0} + \DataV{5}{6}{3}{1} + \DataV{5}{6}{6}{3} \\
&\quad {}+ \DataV{5}{9}{9}{2} + \DataV{5}{3}{2}{0} \\
\DataV{2}{1}{0}{0}\otimes\DataV{3}{0}{1}{0} &= \DataV{2}{1}{0}{0} + \DataV{4}{1}{1}{0} + \DataV{4}{1}{1}{1} + 2\DataV{3}{1}{1}{0} \\
&\quad {}+ 2\DataV{4}{2}{3}{1} + \DataV{5}{6}{12}{6} + 3\DataV{3}{0}{1}{0} + 2\DataV{4}{0}{1}{0} \\
&\quad {}+ \DataV{5}{1}{6}{3} + \DataV{2}{0}{1}{0} + 2\DataV{4}{1}{2}{0} + 2\DataV{4}{1}{3}{2} \\
&\quad {}+ 2\DataV{5}{9}{15}{6} + \DataV{5}{3}{11}{6} + \DataV{4}{0}{1}{1} + \DataV{4}{2}{1}{0} \\
&\quad {}+ 2\DataV{3}{1}{1}{1} + \DataV{5}{6}{11}{3} + 2\DataV{4}{3}{2}{1} + 3\DataV{4}{2}{4}{2} \\
&\quad {}+ \DataV{5}{6}{15}{9} + \DataV{4}{0}{2}{1} + 2\DataV{3}{0}{1}{1} + \DataV{5}{3}{9}{3} \\
&\quad {}+ \DataV{2}{0}{0}{1} + \DataV{5}{8}{9}{3} + \DataV{4}{1}{2}{3} + \DataV{5}{1}{3}{0} \\
&\quad {}+ \DataV{5}{3}{6}{1} + \DataV{5}{9}{9}{2} + \DataV{5}{3}{3}{0} \\
\DataV{2}{1}{0}{0}\otimes\DataV{3}{0}{0}{1} &= \DataV{3}{1}{1}{1} + \DataV{4}{2}{4}{2} + \DataV{4}{0}{2}{1} + \DataV{3}{0}{1}{1} \\
&\quad {}+ \DataV{4}{1}{2}{3} + \DataV{5}{6}{8}{6} + \DataV{5}{3}{9}{8} + \DataV{4}{0}{1}{2} \\
&\quad {}+ \DataV{3}{0}{0}{1} + \DataV{5}{1}{3}{6} \\
\DataV{2}{1}{0}{0}\otimes\DataV{3}{0}{1}{1} &= \DataV{3}{1}{1}{0} + \DataV{4}{2}{3}{1} + 2\DataV{3}{0}{1}{0} + \DataV{4}{0}{1}{0} \\
&\quad {}+ \DataV{2}{0}{1}{0} + \DataV{4}{1}{2}{0} + 2\DataV{4}{1}{3}{2} + \DataV{5}{9}{15}{6} \\
&\quad {}+ \DataV{5}{3}{11}{6} + \DataV{4}{0}{1}{1} + 2\DataV{3}{1}{1}{1} + \DataV{4}{3}{2}{1} \\
&\quad {}+ 3\DataV{4}{2}{4}{2} + \DataV{5}{6}{15}{9} + 2\DataV{4}{0}{2}{1} + 3\DataV{3}{0}{1}{1} \\
&\quad {}+ \DataV{5}{3}{9}{3} + \DataV{5}{2}{9}{9} + \DataV{2}{0}{0}{1} + \DataV{5}{8}{9}{3} \\
&\quad {}+ 2\DataV{4}{1}{2}{3} + \DataV{5}{6}{8}{6} + \DataV{5}{3}{9}{8} + \DataV{4}{0}{1}{2} \\
&\quad {}+ \DataV{3}{0}{0}{1} \\
\DataV{2}{1}{0}{0}\otimes\DataV{3}{1}{1}{0} &= \DataV{4}{1}{0}{0} + 2\DataV{2}{1}{0}{0} + 2\DataV{4}{1}{1}{0} + \DataV{4}{1}{1}{1} \\
&\quad {}+ 2\DataV{3}{1}{0}{0} + 4\DataV{3}{1}{1}{0} + 3\DataV{4}{2}{3}{1} + \DataV{5}{6}{12}{6} \\
&\quad {}+ 2\DataV{3}{0}{1}{0} + \DataV{4}{0}{1}{0} + 2\DataV{2}{0}{1}{0} + 3\DataV{4}{1}{2}{0} \\
&\quad {}+ \DataV{4}{1}{3}{2} + \DataV{5}{9}{15}{6} + L + 3\DataV{4}{2}{1}{0} \\
&\quad {}+ 3\DataV{3}{1}{1}{1} + \DataV{5}{6}{11}{3} + 3\DataV{4}{3}{2}{1} + 2\DataV{4}{2}{4}{2} \\
&\quad {}+ \DataV{3}{0}{1}{1} + \DataV{2}{0}{0}{1} + \DataV{5}{8}{9}{3} + \DataV{5}{3}{1}{0} \\
&\quad {}+ \DataV{5}{6}{3}{1} + \DataV{5}{3}{6}{1} + \DataV{5}{6}{6}{3} + 2\DataV{5}{9}{9}{2} \\
&\quad {}+ \DataV{5}{3}{2}{0} + \DataV{5}{2}{3}{0} + \DataV{5}{3}{3}{0} \\
\DataV{2}{1}{0}{0}\otimes\DataV{3}{1}{1}{1} &= \DataV{2}{1}{0}{0} + \DataV{4}{1}{1}{0} + \DataV{4}{1}{1}{1} + \DataV{3}{1}{0}{0} \\
&\quad {}+ 3\DataV{3}{1}{1}{0} + 3\DataV{4}{2}{3}{1} + \DataV{5}{6}{12}{6} + 2\DataV{3}{0}{1}{0} \\
&\quad {}+ \DataV{4}{0}{1}{0} + 2\DataV{2}{0}{1}{0} + 2\DataV{4}{1}{2}{0} + 2\DataV{4}{1}{3}{2} \\
&\quad {}+ \DataV{5}{9}{15}{6} + L + \DataV{4}{2}{1}{0} + 4\DataV{3}{1}{1}{1} \\
&\quad {}+ \DataV{5}{6}{11}{3} + \DataV{5}{3}{6}{6} + 2\DataV{4}{3}{2}{1} + 3\DataV{4}{2}{4}{2} \\
&\quad {}+ \DataV{5}{6}{15}{9} + \DataV{4}{0}{2}{1} + 2\DataV{3}{0}{1}{1} + \DataV{2}{0}{0}{1} \\
&\quad {}+ \DataV{5}{8}{9}{3} + \DataV{4}{1}{2}{3} + \DataV{5}{6}{8}{6} + \DataV{3}{0}{0}{1} \\
&\quad {}+ \DataV{5}{6}{3}{1} + \DataV{5}{6}{6}{3} + \DataV{5}{9}{9}{2} \\
\DataV{2}{0}{1}{0}\otimes\DataV{3}{1}{0}{0} &= \DataV{4}{1}{0}{0} + \DataV{2}{1}{0}{0} + \DataV{4}{1}{1}{0} + \DataV{3}{1}{0}{0} \\
&\quad {}+ 2\DataV{3}{1}{1}{0} + \DataV{4}{2}{3}{1} + \DataV{3}{0}{1}{0} + \DataV{2}{0}{1}{0} \\
&\quad {}+ \DataV{4}{1}{2}{0} + \DataV{5}{9}{15}{6} + 2\DataV{4}{2}{1}{0} + \DataV{3}{1}{1}{1} \\
&\quad {}+ 2\DataV{4}{3}{2}{1} + \DataV{4}{2}{4}{2} + \DataV{5}{8}{9}{3} + \DataV{5}{3}{1}{0} \\
&\quad {}+ \DataV{5}{6}{3}{1} + \DataV{5}{9}{9}{2} + \DataV{5}{3}{3}{0} \\
\DataV{2}{0}{1}{0}\otimes\DataV{3}{0}{1}{0} &= \DataV{2}{1}{0}{0} + \DataV{4}{1}{1}{0} + \DataV{4}{1}{1}{1} + \DataV{3}{1}{0}{0} \\
&\quad {}+ 3\DataV{3}{1}{1}{0} + 4\DataV{4}{2}{3}{1} + 2\DataV{5}{6}{12}{6} + 3\DataV{3}{0}{1}{0} \\
&\quad {}+ 3\DataV{4}{0}{1}{0} + \DataV{5}{1}{6}{3} + 2\DataV{2}{0}{1}{0} + 3\DataV{4}{1}{2}{0} \\
&\quad {}+ 4\DataV{4}{1}{3}{2} + 2\DataV{5}{9}{15}{6} + 2\DataV{5}{3}{11}{6} + \DataV{5}{0}{3}{2} \\
&\quad {}+ \DataV{4}{0}{1}{1} + L + \DataV{4}{2}{1}{0} + 4\DataV{3}{1}{1}{1} \\
&\quad {}+ 2\DataV{5}{6}{11}{3} + \DataV{5}{3}{6}{6} + 2\DataV{4}{3}{2}{1} + 5\DataV{4}{2}{4}{2} \\
&\quad {}+ 2\DataV{5}{6}{15}{9} + 3\DataV{4}{0}{2}{1} + \DataV{5}{0}{3}{1} + 3\DataV{3}{0}{1}{1} \\
&\quad {}+ 2\DataV{5}{3}{9}{3} + \DataV{5}{2}{9}{9} + \DataV{2}{0}{0}{1} + \DataV{5}{8}{9}{3} \\
&\quad {}+ 2\DataV{4}{1}{2}{3} + \DataV{5}{6}{8}{6} + \DataV{5}{3}{9}{8} + \DataV{4}{0}{1}{2} \\
&\quad {}+ \DataV{3}{0}{0}{1} + \DataV{5}{0}{1}{0} + \DataV{5}{1}{3}{0} + \DataV{5}{3}{6}{1} \\
&\quad {}+ \DataV{5}{6}{6}{3} + \DataV{5}{9}{9}{2} + \DataV{5}{2}{3}{0} \\
\DataV{2}{0}{1}{0}\otimes\DataV{3}{0}{0}{1} &= \DataV{3}{0}{1}{0} + \DataV{2}{0}{1}{0} + \DataV{4}{1}{3}{2} + \DataV{4}{0}{1}{1} \\
&\quad {}+ \DataV{3}{1}{1}{1} + \DataV{4}{2}{4}{2} + \DataV{5}{6}{15}{9} + \DataV{4}{0}{2}{1} \\
&\quad {}+ 2\DataV{3}{0}{1}{1} + \DataV{5}{2}{9}{9} + \DataV{2}{0}{0}{1} + 2\DataV{4}{1}{2}{3} \\
&\quad {}+ \DataV{5}{3}{9}{8} + \DataV{5}{0}{3}{3} + 2\DataV{4}{0}{1}{2} + \DataV{3}{0}{0}{1} \\
&\quad {}+ \DataV{5}{1}{3}{6} + \DataV{5}{0}{1}{3} + \DataV{4}{0}{0}{1} \\
\DataV{2}{0}{1}{0}\otimes\DataV{3}{0}{1}{1} &= \DataV{2}{1}{0}{0} + \DataV{4}{1}{1}{1} + 2\DataV{3}{1}{1}{0} + 2\DataV{4}{2}{3}{1} \\
&\quad {}+ \DataV{5}{6}{12}{6} + 3\DataV{3}{0}{1}{0} + 2\DataV{4}{0}{1}{0} + \DataV{5}{1}{6}{3} \\
&\quad {}+ 2\DataV{2}{0}{1}{0} + \DataV{4}{1}{2}{0} + 4\DataV{4}{1}{3}{2} + \DataV{5}{9}{15}{6} \\
&\quad {}+ 2\DataV{5}{3}{11}{6} + \DataV{5}{0}{3}{2} + 2\DataV{4}{0}{1}{1} + L \\
&\quad {}+ 4\DataV{3}{1}{1}{1} + \DataV{5}{6}{11}{3} + \DataV{5}{3}{6}{6} + \DataV{4}{3}{2}{1} \\
&\quad {}+ 4\DataV{4}{2}{4}{2} + 2\DataV{5}{6}{15}{9} + 4\DataV{4}{0}{2}{1} + \DataV{5}{0}{3}{1} \\
&\quad {}+ 5\DataV{3}{0}{1}{1} + \DataV{5}{3}{9}{3} + 2\DataV{5}{2}{9}{9} + \DataV{5}{0}{2}{3} \\
&\quad {}+ 2\DataV{2}{0}{0}{1} + 4\DataV{4}{1}{2}{3} + \DataV{5}{6}{8}{6} + 2\DataV{5}{3}{9}{8} \\
&\quad {}+ \DataV{5}{0}{3}{3} + 3\DataV{4}{0}{1}{2} + 2\DataV{3}{0}{0}{1} + \DataV{5}{1}{3}{6} \\
&\quad {}+ \DataV{4}{0}{0}{1} \\
\DataV{2}{0}{1}{0}\otimes\DataV{3}{1}{1}{0} &= \DataV{4}{1}{0}{0} + 2\DataV{2}{1}{0}{0} + 2\DataV{4}{1}{1}{0} + \DataV{4}{1}{1}{1} \\
&\quad {}+ 2\DataV{3}{1}{0}{0} + 5\DataV{3}{1}{1}{0} + 4\DataV{4}{2}{3}{1} + \DataV{5}{6}{12}{6} \\
&\quad {}+ 3\DataV{3}{0}{1}{0} + 2\DataV{4}{0}{1}{0} + 2\DataV{2}{0}{1}{0} + 4\DataV{4}{1}{2}{0} \\
&\quad {}+ 2\DataV{4}{1}{3}{2} + 2\DataV{5}{9}{15}{6} + \DataV{5}{3}{11}{6} + L \\
&\quad {}+ 3\DataV{4}{2}{1}{0} + 4\DataV{3}{1}{1}{1} + 2\DataV{5}{6}{11}{3} + 4\DataV{4}{3}{2}{1} \\
&\quad {}+ 4\DataV{4}{2}{4}{2} + \DataV{5}{6}{15}{9} + \DataV{4}{0}{2}{1} + 2\DataV{3}{0}{1}{1} \\
&\quad {}+ \DataV{5}{3}{9}{3} + \DataV{2}{0}{0}{1} + 2\DataV{5}{8}{9}{3} + \DataV{4}{1}{2}{3} \\
&\quad {}+ \DataV{5}{6}{8}{6} + \DataV{5}{6}{3}{1} + \DataV{5}{1}{3}{0} + \DataV{5}{3}{6}{1} \\
&\quad {}+ \DataV{5}{6}{6}{3} + 2\DataV{5}{9}{9}{2} + \DataV{5}{3}{2}{0} + \DataV{5}{2}{3}{0} \\
&\quad {}+ \DataV{5}{3}{3}{0} \\
\DataV{2}{0}{1}{0}\otimes\DataV{3}{1}{1}{1} &= 2\DataV{2}{1}{0}{0} + \DataV{4}{1}{1}{0} + 2\DataV{4}{1}{1}{1} + \DataV{3}{1}{0}{0} \\
&\quad {}+ 4\DataV{3}{1}{1}{0} + 4\DataV{4}{2}{3}{1} + 2\DataV{5}{6}{12}{6} + 4\DataV{3}{0}{1}{0} \\
&\quad {}+ 2\DataV{4}{0}{1}{0} + \DataV{5}{1}{6}{3} + 3\DataV{2}{0}{1}{0} + 2\DataV{4}{1}{2}{0} \\
&\quad {}+ 4\DataV{4}{1}{3}{2} + 2\DataV{5}{9}{15}{6} + \DataV{5}{3}{11}{6} + \DataV{4}{0}{1}{1} \\
&\quad {}+ L + \DataV{4}{2}{1}{0} + 5\DataV{3}{1}{1}{1} + \DataV{5}{6}{11}{3} \\
&\quad {}+ \DataV{5}{3}{6}{6} + 3\DataV{4}{3}{2}{1} + 5\DataV{4}{2}{4}{2} + 2\DataV{5}{6}{15}{9} \\
&\quad {}+ 2\DataV{4}{0}{2}{1} + 4\DataV{3}{0}{1}{1} + \DataV{5}{3}{9}{3} + \DataV{5}{2}{9}{9} \\
&\quad {}+ 2\DataV{2}{0}{0}{1} + \DataV{5}{8}{9}{3} + 3\DataV{4}{1}{2}{3} + \DataV{5}{6}{8}{6} \\
&\quad {}+ \DataV{5}{3}{9}{8} + \DataV{4}{0}{1}{2} + \DataV{3}{0}{0}{1} + \DataV{5}{3}{6}{1} \\
&\quad {}+ \DataV{5}{6}{6}{3} + \DataV{5}{9}{9}{2} \\
\DataV{2}{0}{0}{1}\otimes\DataV{3}{1}{0}{0} &= \DataV{3}{1}{0}{0} + \DataV{3}{1}{1}{0} + \DataV{4}{1}{2}{0} + \DataV{4}{2}{1}{0} \\
&\quad {}+ \DataV{3}{1}{1}{1} + \DataV{4}{3}{2}{1} + \DataV{4}{2}{4}{2} + \DataV{5}{8}{9}{3} \\
&\quad {}+ \DataV{5}{6}{8}{6} + \DataV{5}{6}{3}{1} \\
\DataV{2}{0}{0}{1}\otimes\DataV{3}{0}{1}{0} &= \DataV{2}{1}{0}{0} + \DataV{4}{1}{1}{0} + \DataV{4}{1}{1}{1} + 2\DataV{3}{1}{1}{0} \\
&\quad {}+ 2\DataV{4}{2}{3}{1} + \DataV{5}{6}{12}{6} + 3\DataV{3}{0}{1}{0} + 2\DataV{4}{0}{1}{0} \\
&\quad {}+ \DataV{5}{1}{6}{3} + \DataV{2}{0}{1}{0} + \DataV{4}{1}{2}{0} + 2\DataV{4}{1}{3}{2} \\
&\quad {}+ \DataV{5}{9}{15}{6} + \DataV{5}{3}{11}{6} + \DataV{4}{0}{1}{1} + 2\DataV{3}{1}{1}{1} \\
&\quad {}+ \DataV{5}{6}{11}{3} + \DataV{4}{3}{2}{1} + 3\DataV{4}{2}{4}{2} + 2\DataV{5}{6}{15}{9} \\
&\quad {}+ 2\DataV{4}{0}{2}{1} + \DataV{5}{0}{3}{1} + 2\DataV{3}{0}{1}{1} + \DataV{5}{3}{9}{3} \\
&\quad {}+ \DataV{5}{2}{9}{9} + \DataV{2}{0}{0}{1} + 2\DataV{4}{1}{2}{3} + \DataV{5}{3}{9}{8} \\
&\quad {}+ \DataV{5}{0}{3}{3} + \DataV{4}{0}{1}{2} + \DataV{5}{3}{6}{1} \\
\DataV{2}{0}{0}{1}\otimes\DataV{3}{0}{0}{1} &= \DataV{2}{0}{1}{0} + \DataV{4}{1}{3}{2} + \DataV{4}{0}{1}{1} + L \\
&\quad {}+ \DataV{3}{1}{1}{1} + \DataV{5}{3}{6}{6} + \DataV{4}{0}{2}{1} + 2\DataV{3}{0}{1}{1} \\
&\quad {}+ \DataV{5}{2}{9}{9} + \DataV{5}{0}{2}{3} + \DataV{2}{0}{0}{1} + \DataV{4}{1}{2}{3} \\
&\quad {}+ 2\DataV{4}{0}{1}{2} + 2\DataV{3}{0}{0}{1} + \DataV{5}{1}{3}{6} + \DataV{5}{0}{1}{3} \\
&\quad {}+ \DataV{4}{0}{0}{1} + \DataV{5}{0}{0}{1} \\
\DataV{2}{0}{0}{1}\otimes\DataV{3}{0}{1}{1} &= \DataV{2}{1}{0}{0} + \DataV{4}{1}{1}{1} + \DataV{3}{1}{1}{0} + \DataV{4}{2}{3}{1} \\
&\quad {}+ \DataV{5}{6}{12}{6} + 2\DataV{3}{0}{1}{0} + \DataV{4}{0}{1}{0} + \DataV{5}{1}{6}{3} \\
&\quad {}+ 2\DataV{2}{0}{1}{0} + 3\DataV{4}{1}{3}{2} + \DataV{5}{3}{11}{6} + \DataV{5}{0}{3}{2} \\
&\quad {}+ 2\DataV{4}{0}{1}{1} + L + 3\DataV{3}{1}{1}{1} + \DataV{5}{3}{6}{6} \\
&\quad {}+ 2\DataV{4}{2}{4}{2} + \DataV{5}{6}{15}{9} + 3\DataV{4}{0}{2}{1} + 4\DataV{3}{0}{1}{1} \\
&\quad {}+ 2\DataV{5}{2}{9}{9} + \DataV{5}{0}{2}{3} + 2\DataV{2}{0}{0}{1} + 3\DataV{4}{1}{2}{3} \\
&\quad {}+ \DataV{5}{3}{9}{8} + \DataV{5}{0}{3}{3} + 3\DataV{4}{0}{1}{2} + 2\DataV{3}{0}{0}{1} \\
&\quad {}+ \DataV{5}{1}{3}{6} + \DataV{5}{0}{1}{3} + \DataV{4}{0}{0}{1} \\
\DataV{2}{0}{0}{1}\otimes\DataV{3}{1}{1}{0} &= \DataV{2}{1}{0}{0} + \DataV{4}{1}{1}{0} + \DataV{3}{1}{0}{0} + 3\DataV{3}{1}{1}{0} \\
&\quad {}+ 2\DataV{4}{2}{3}{1} + 2\DataV{3}{0}{1}{0} + \DataV{4}{0}{1}{0} + \DataV{2}{0}{1}{0} \\
&\quad {}+ 2\DataV{4}{1}{2}{0} + \DataV{4}{1}{3}{2} + \DataV{5}{9}{15}{6} + \DataV{4}{2}{1}{0} \\
&\quad {}+ 2\DataV{3}{1}{1}{1} + \DataV{5}{6}{11}{3} + 2\DataV{4}{3}{2}{1} + 3\DataV{4}{2}{4}{2} \\
&\quad {}+ \DataV{5}{6}{15}{9} + \DataV{4}{0}{2}{1} + \DataV{3}{0}{1}{1} + \DataV{5}{3}{9}{3} \\
&\quad {}+ \DataV{5}{8}{9}{3} + \DataV{4}{1}{2}{3} + \DataV{5}{6}{8}{6} + \DataV{5}{3}{9}{8} \\
&\quad {}+ \DataV{5}{9}{9}{2} \\
\DataV{2}{0}{0}{1}\otimes\DataV{3}{1}{1}{1} &= \DataV{2}{1}{0}{0} + \DataV{4}{1}{1}{1} + \DataV{3}{1}{0}{0} + 2\DataV{3}{1}{1}{0} \\
&\quad {}+ 2\DataV{4}{2}{3}{1} + \DataV{5}{6}{12}{6} + 2\DataV{3}{0}{1}{0} + \DataV{4}{0}{1}{0} \\
&\quad {}+ 2\DataV{2}{0}{1}{0} + \DataV{4}{1}{2}{0} + 3\DataV{4}{1}{3}{2} + \DataV{5}{9}{15}{6} \\
&\quad {}+ \DataV{5}{3}{11}{6} + \DataV{4}{0}{1}{1} + L + 4\DataV{3}{1}{1}{1} \\
&\quad {}+ \DataV{5}{3}{6}{6} + \DataV{4}{3}{2}{1} + 3\DataV{4}{2}{4}{2} + \DataV{5}{6}{15}{9} \\
&\quad {}+ 2\DataV{4}{0}{2}{1} + 3\DataV{3}{0}{1}{1} + \DataV{5}{2}{9}{9} + \DataV{2}{0}{0}{1} \\
&\quad {}+ 2\DataV{4}{1}{2}{3} + \DataV{5}{6}{8}{6} + \DataV{5}{3}{9}{8} + \DataV{4}{0}{1}{2} \\
&\quad {}+ \DataV{3}{0}{0}{1} + \DataV{5}{1}{3}{6} + \DataV{5}{6}{6}{3} \\
\end{align*}

%% file: qdim/level0.tex
\scriptsize
\begin{longtable}{@{}>{$}l<{$} >{$\displaystyle}p{0.39\linewidth}<{$} >{$\displaystyle}p{0.39\linewidth}<{$}@{}}
\toprule
V & N(V) & D(V) \\
\midrule\endhead
\DataUnit & 1 & 1 \\
\bottomrule
\end{longtable}
\normalsize

%% file: qdim/level1.tex
\scriptsize
\begin{longtable}{@{}>{$}l<{$} >{$\displaystyle}p{0.39\linewidth}<{$} >{$\displaystyle}p{0.39\linewidth}<{$}@{}}
\toprule
V & N(V) & D(V) \\
\midrule\endhead
L & \left[6\nv + 1\right]\allowbreak\, \left[5\nv - 1\right]\allowbreak\, \left[4\nv\right] & \left[\nv + 1\right]\allowbreak\, \left[2\nv\right]\, \left[1\right] \\
\bottomrule
\end{longtable}
\normalsize

%% file: qdim/level2.tex
\scriptsize
\begin{longtable}{@{}>{$}l<{$} >{$\displaystyle}p{0.39\linewidth}<{$} >{$\displaystyle}p{0.39\linewidth}<{$}@{}}
\toprule
V & N(V) & D(V) \\
\midrule\endhead
\DataV{2}{1}{0}{0} & \left[6\nv + 1\right]\allowbreak\, \left[6\nv\right]\allowbreak\, \left[5\nv\right]\allowbreak\, \left[4\nv\right]\allowbreak\, \left[3\nv - 1\right]\allowbreak\, \left[3\nv - 3\right] & \left[2\nv + 2\right]\allowbreak\, \left[2\nv\right]\allowbreak\, \left[\nv + 2\right]\allowbreak\, \left[\nv + 1\right]\allowbreak\, \left[\nv - 1\right]\allowbreak\, \left[1\right] \\
\DataV{2}{0}{1}{0} & \left[3\nv - 1\right]\allowbreak\, \left[4\nv - 2\right]\allowbreak\, \left[5\nv\right]\allowbreak\, \left[4\nv + 1\right]\allowbreak\, \left[5\nv - 1\right]\allowbreak\, \left[6\nv + 1\right]\allowbreak\, \left[6\nv + 2\right] & \left[1\right]\allowbreak\, \left[2\right]\allowbreak\, \left[\nv\right]\allowbreak\, \left[2\nv - 1\right]\allowbreak\, \left[\nv + 1\right]\allowbreak\, \left[2\nv + 2\right]\allowbreak\, \left[3\nv + 1\right] \\
\DataV{2}{0}{0}{1} & \left[6\nv + 3\right]\allowbreak\, \left[6\nv\right]\allowbreak\, \left[5\nv\right]\allowbreak\, \left[5\nv - 1\right]\allowbreak\, \left[4\nv + 1\right]\allowbreak\, \left[4\nv\right] & \left[2\nv + 1\right]\allowbreak\, \left[2\nv\right]\allowbreak\, \left[\nv + 2\right]\allowbreak\, \left[\nv + 1\right]\allowbreak\, \left[2\right]\allowbreak\, \left[1\right] \\
\bottomrule
\end{longtable}
\normalsize

%% file: qdim/level3.tex
\scriptsize
\begin{longtable}{@{}>{$}l<{$} >{$\displaystyle}p{0.39\linewidth}<{$} >{$\displaystyle}p{0.39\linewidth}<{$}@{}}
\toprule
V & N(V) & D(V) \\
\midrule\endhead
\DataV{3}{1}{0}{0} & -\left[\nv - 5\right]\allowbreak\, \left[5\nv - 1\right]\allowbreak\, \left[6\nv\right]\allowbreak\, \left[4\nv - 1\right]\allowbreak\, \left[5\nv\right]\allowbreak\, \left[6\nv + 1\right]\allowbreak\, \left[2\nv - 2\right]\allowbreak\, \left[3\nv - 1\right]\allowbreak\, \left[4\nv\right] & \left[2\right]\allowbreak\, \left[\nv - 1\right]\allowbreak\, \left[2\nv\right]\allowbreak\, \left[2\nv + 3\right]\allowbreak\, \left[\nv + 2\right]\allowbreak\, \left[1\right]\allowbreak\, \left[3\nv + 3\right]\allowbreak\, \left[2\nv + 2\right]\allowbreak\, \left[\nv + 1\right] \\
\DataV{3}{0}{1}{0} & \left[2\nv - 1\right]\allowbreak\, \left[3\nv - 3\right]\allowbreak\, \left[3\nv + 1\right]\allowbreak\, \left[4\nv - 1\right]\allowbreak\, \left[4\nv\right]\allowbreak\, \left[4\nv - 2\right]\allowbreak\, \left[5\nv\right]\allowbreak\, \left[5\nv - 1\right]\allowbreak\, \left[5\nv + 1\right]\allowbreak\, \left[6\nv + 1\right]\allowbreak\, \left[6\nv + 2\right]\allowbreak\, \left[6\nv + 3\right] & \left[1\right]\allowbreak\, \left[2\right]\allowbreak\, \left[3\right]\allowbreak\, \left[\nv - 1\right]\allowbreak\, \left[\nv\right]\allowbreak\, \left[\nv + 1\right]\allowbreak\, \left[2\nv - 1\right]\allowbreak\, \left[2\nv\right]\allowbreak\, \left[2\nv + 1\right]\allowbreak\, \left[2\nv + 2\right]\allowbreak\, \left[3\nv + 1\right]\allowbreak\, \left[3\nv + 3\right] \\
\DataV{3}{0}{0}{1} & \left[6\nv + 5\right]\allowbreak\, \left[6\nv + 1\right]\allowbreak\, \left[6\nv\right]\allowbreak\, \left[5\nv + 1\right]\allowbreak\, \left[5\nv\right]\allowbreak\, \left[5\nv - 1\right]\allowbreak\, \left[4\nv + 2\right]\allowbreak\, \left[4\nv + 1\right]\allowbreak\, \left[4\nv\right] & \left[2\nv + 2\right]\allowbreak\, \left[2\nv + 1\right]\allowbreak\, \left[2\nv\right]\allowbreak\, \left[\nv + 3\right]\allowbreak\, \left[\nv + 2\right]\allowbreak\, \left[\nv + 1\right]\allowbreak\, \left[3\right]\allowbreak\, \left[2\right]\allowbreak\, \left[1\right] \\
\DataV{3}{0}{1}{1} & \left[3\nv - 1\right]\allowbreak\, \left[4\nv - 2\right]\allowbreak\, \left[4\nv\right]\allowbreak\, \left[4\nv + 2\right]\allowbreak\, \left[5\nv - 1\right]\allowbreak\, \left[5\nv\right]\allowbreak\, \left[5\nv + 1\right]\allowbreak\, \left[6\nv\right]\allowbreak\, \left[6\nv + 2\right]\allowbreak\, \left[6\nv + 4\right] & \left[1\right]\allowbreak\, \left[1\right]\allowbreak\, \left[3\right]\allowbreak\, \left[\nv\right]\allowbreak\, \left[\nv + 1\right]\allowbreak\, \left[\nv + 2\right]\allowbreak\, \left[2\nv - 1\right]\allowbreak\, \left[2\nv + 1\right]\allowbreak\, \left[2\nv + 3\right]\allowbreak\, \left[3\nv + 2\right] \\
\DataV{3}{1}{1}{0} & \left[2\nv - 4\right]\allowbreak\, \left[2\nv - 2\right]\allowbreak\, \left[4\nv - 2\right]\allowbreak\, \left[4\nv - 1\right]\allowbreak\, \left[4\nv\right]\allowbreak\, \left[4\nv + 1\right]\allowbreak\, \left[5\nv\right]\allowbreak\, \left[6\nv\right]\allowbreak\, \left[6\nv + 1\right]\allowbreak\, \left[6\nv + 2\right] & \left[1\right]\allowbreak\, \left[\nv - 2\right]\allowbreak\, \left[\nv - 1\right]\allowbreak\, \left[\nv\right]\allowbreak\, \left[\nv + 1\right]\allowbreak\, \left[\nv + 1\right]\allowbreak\, \left[\nv + 2\right]\allowbreak\, \left[\nv + 3\right]\allowbreak\, \left[3\nv + 1\right]\allowbreak\, \left[3\nv + 3\right] \\
\DataV{3}{1}{1}{1} & \left[6\nv\right]\allowbreak\, \left[6\nv + 1\right]\allowbreak\, \left[6\nv + 3\right]\allowbreak\, \left[5\nv + 1\right]\allowbreak\, \left[4\nv - 1\right]\allowbreak\, \left[4\nv + 1\right]\allowbreak\, \left[5\nv - 1\right]\allowbreak\, \left[3\nv - 1\right]\allowbreak\, \left[3\nv - 3\right] & \left[\nv + 1\right]\allowbreak\, \left[\nv + 1\right]\allowbreak\, \left[2\nv + 3\right]\allowbreak\, \left[2\nv + 1\right]\allowbreak\, \left[2\nv\right]\allowbreak\, \left[\nv - 1\right]\allowbreak\, \left[\nv + 3\right] \\
\bottomrule
\end{longtable}
\normalsize

%% file: qdim/level4.tex
\scriptsize
\begin{longtable}{@{}>{$}l<{$} >{$\displaystyle}p{0.39\linewidth}<{$} >{$\displaystyle}p{0.39\linewidth}<{$}@{}}
\toprule
V & N(V) & D(V) \\
\midrule\endhead
\DataV{4}{1}{0}{0} & -\left[6\nv + 1\right]\allowbreak\, \left[6\nv\right]\allowbreak\, \left[5\nv\right]\allowbreak\, \left[5\nv - 1\right]\allowbreak\, \left[4\nv\right]\allowbreak\, \left[4\nv - 1\right]\allowbreak\, \left[4\nv - 2\right]\allowbreak\, \left[3\nv - 1\right]\allowbreak\, \left[3\nv - 2\right]\allowbreak\, \left[2\nv - 2\right]\allowbreak\, \left[\nv + 7\right]\allowbreak\, \left[\nv - 3\right] & \left[4\nv + 4\right]\allowbreak\, \left[3\nv + 4\right]\allowbreak\, \left[3\nv + 3\right]\allowbreak\, \left[2\nv + 3\right]\allowbreak\, \left[2\nv + 2\right]\allowbreak\, \left[2\nv\right]\allowbreak\, \left[\nv + 3\right]\allowbreak\, \left[\nv + 2\right]\allowbreak\, \left[\nv + 1\right]\allowbreak\, \left[\nv - 1\right]\allowbreak\, \left[2\right]\allowbreak\, \left[1\right] \\
\DataV{4}{0}{0}{1} & \left[6\nv + 7\right]\allowbreak\, \left[6\nv + 2\right]\allowbreak\, \left[6\nv + 1\right]\allowbreak\, \left[6\nv\right]\allowbreak\, \left[5\nv + 2\right]\allowbreak\, \left[5\nv + 1\right]\allowbreak\, \left[5\nv\right]\allowbreak\, \left[5\nv - 1\right]\allowbreak\, \left[4\nv + 3\right]\allowbreak\, \left[4\nv + 2\right]\allowbreak\, \left[4\nv + 1\right]\allowbreak\, \left[4\nv\right] & \left[2\nv + 3\right]\allowbreak\, \left[2\nv + 2\right]\allowbreak\, \left[2\nv + 1\right]\allowbreak\, \left[2\nv\right]\allowbreak\, \left[\nv + 4\right]\allowbreak\, \left[\nv + 3\right]\allowbreak\, \left[\nv + 2\right]\allowbreak\, \left[\nv + 1\right]\allowbreak\, \left[4\right]\allowbreak\, \left[3\right]\allowbreak\, \left[2\right]\allowbreak\, \left[1\right] \\
\DataV{4}{0}{1}{0} & \left[3\nv - 3\right]\allowbreak\, \left[3\nv - 1\right]\allowbreak\, \left[3\nv - 2\right]\allowbreak\, \left[2\nv - 4\right]\allowbreak\, \left[4\nv - 1\right]\allowbreak\, \left[4\nv - 2\right]\allowbreak\, \left[5\nv\right]\allowbreak\, \left[4\nv\right]\allowbreak\, \left[4\nv\right]\allowbreak\, \left[4\nv + 1\right]\allowbreak\, \left[5\nv - 1\right]\allowbreak\, \left[6\nv + 3\right]\allowbreak\, \left[5\nv + 1\right]\allowbreak\, \left[5\nv + 2\right]\allowbreak\, \left[6\nv + 1\right]\allowbreak\, \left[6\nv + 2\right]\allowbreak\, \left[6\nv + 4\right] & \left[1\right]\allowbreak\, \left[2\right]\allowbreak\, \left[3\right]\allowbreak\, \left[4\right]\allowbreak\, \left[\nv + 1\right]\allowbreak\, \left[2\nv + 2\right]\allowbreak\, \left[3\nv + 3\right]\allowbreak\, \left[4\nv + 4\right]\allowbreak\, \left[\nv - 1\right]\allowbreak\, \left[\nv - 2\right]\allowbreak\, \left[2\nv - 1\right]\allowbreak\, \left[\nv\right]\allowbreak\, \left[2\nv\right]\allowbreak\, \left[2\nv\right]\allowbreak\, \left[2\nv + 1\right]\allowbreak\, \left[3\nv + 1\right]\allowbreak\, \left[3\nv + 2\right] \\
\DataV{4}{0}{1}{1} & \left[6\nv + 5\right]\allowbreak\, \left[6\nv + 4\right]\allowbreak\, \left[6\nv + 1\right]\allowbreak\, \left[6\nv\right]\allowbreak\, \left[5\nv + 2\right]\allowbreak\, \left[5\nv + 1\right]\allowbreak\, \left[5\nv\right]\allowbreak\, \left[5\nv - 1\right]\allowbreak\, \left[4\nv + 3\right]\allowbreak\, \left[4\nv\right]\allowbreak\, \left[4\nv - 1\right]\allowbreak\, \left[4\nv - 2\right]\allowbreak\, \left[3\nv\right]\allowbreak\, \left[3\nv - 1\right] & \left[3\nv + 3\right]\allowbreak\, \left[3\nv + 2\right]\allowbreak\, \left[2\nv + 3\right]\allowbreak\, \left[2\nv + 4\right]\allowbreak\, \left[2\nv\right]\allowbreak\, \left[2\nv - 1\right]\allowbreak\, \left[\nv + 2\right]\allowbreak\, \left[\nv + 1\right]\allowbreak\, \left[\nv + 1\right]\allowbreak\, \left[\nv\right]\allowbreak\, \left[3\right]\allowbreak\, \left[2\right]\allowbreak\, \left[2\right]\allowbreak\, \left[1\right] \\
\DataV{4}{1}{1}{0} & \left[\nv - 5\right]\allowbreak\, \left[2\nv - 4\right]\allowbreak\, \left[5\nv - 1\right]\allowbreak\, \left[6\nv\right]\allowbreak\, \left[3\nv - 2\right]\allowbreak\, \left[4\nv - 1\right]\allowbreak\, \left[5\nv\right]\allowbreak\, \left[6\nv + 1\right]\allowbreak\, \left[\nv - 3\right]\allowbreak\, \left[4\nv\right]\allowbreak\, \left[5\nv + 1\right]\allowbreak\, \left[6\nv + 2\right]\allowbreak\, \left[3\nv\right]\allowbreak\, \left[4\nv + 1\right] & \left[3\right]\allowbreak\, \left[\nv - 2\right]\allowbreak\, \left[\nv + 3\right]\allowbreak\, \left[2\nv + 4\right]\allowbreak\, \left[2\nv\right]\allowbreak\, \left[3\nv + 1\right]\allowbreak\, \left[\nv + 2\right]\allowbreak\, \left[1\right]\allowbreak\, \left[1\right]\allowbreak\, \left[\nv\right]\allowbreak\, \left[3\nv + 3\right]\allowbreak\, \left[2\nv + 2\right]\allowbreak\, \left[2\nv + 2\right]\allowbreak\, \left[\nv + 1\right] \\
\DataV{4}{1}{1}{1} & \left[6\nv + 3\right]\allowbreak\, \left[6\nv + 2\right]\allowbreak\, \left[6\nv + 1\right]\allowbreak\, \left[6\nv\right]\allowbreak\, \left[5\nv + 2\right]\allowbreak\, \left[5\nv - 1\right]\allowbreak\, \left[4\nv + 1\right]\allowbreak\, \left[4\nv - 2\right]\allowbreak\, \left[3\nv\right]\allowbreak\, \left[3\nv - 1\right]\allowbreak\, \left[3\nv - 2\right]\allowbreak\, \left[3\nv - 3\right] & \left[2\nv + 3\right]\allowbreak\, \left[2\nv + 2\right]\allowbreak\, \left[2\nv + 1\right]\allowbreak\, \left[2\nv\right]\allowbreak\, \left[\nv + 3\right]\allowbreak\, \left[\nv + 2\right]\allowbreak\, \left[\nv + 2\right]\allowbreak\, \left[\nv + 1\right]\allowbreak\, \left[\nv\right]\allowbreak\, \left[\nv - 1\right]\allowbreak\, \left[2\right]\allowbreak\, \left[1\right] \\
\DataV{4}{1}{2}{0} & -\left[6\nv + 3\right]\allowbreak\, \left[6\nv + 2\right]\allowbreak\, \left[6\nv + 1\right]\allowbreak\, \left[6\nv\right]\allowbreak\, \left[5\nv + 1\right]\allowbreak\, \left[5\nv\right]\allowbreak\, \left[4\nv\right]\allowbreak\, \left[4\nv - 1\right]\allowbreak\, \left[4\nv - 2\right]\allowbreak\, \left[3\nv\right]\allowbreak\, \left[3\nv - 1\right]\allowbreak\, \left[3\nv - 2\right]\allowbreak\, \left[3\nv - 3\right]\allowbreak\, \left[\nv - 5\right] & \left[4\nv + 4\right]\allowbreak\, \left[2\nv + 2\right]\allowbreak\, \left[2\nv + 1\right]\allowbreak\, \left[2\nv\right]\allowbreak\, \left[2\nv - 1\right]\allowbreak\, \left[\nv + 4\right]\allowbreak\, \left[\nv + 3\right]\allowbreak\, \left[\nv + 2\right]\allowbreak\, \left[\nv + 1\right]\allowbreak\, \left[\nv + 1\right]\allowbreak\, \left[\nv\right]\allowbreak\, \left[\nv - 1\right]\allowbreak\, \left[2\right]\allowbreak\, \left[1\right] \\
\DataV{4}{2}{1}{0} & -\left[6\nv + 2\right]\allowbreak\, \left[6\nv + 1\right]\allowbreak\, \left[6\nv\right]\allowbreak\, \left[5\nv\right]\allowbreak\, \left[5\nv - 1\right]\allowbreak\, \left[4\nv + 1\right]\allowbreak\, \left[4\nv\right]\allowbreak\, \left[4\nv - 1\right]\allowbreak\, \left[3\nv - 1\right]\allowbreak\, \left[3\nv - 2\right]\allowbreak\, \left[3\nv - 3\right]\allowbreak\, \left[\nv - 3\right]\allowbreak\, \left[6\right] & \left[4\nv + 4\right]\allowbreak\, \left[2\nv + 4\right]\allowbreak\, \left[3\nv + 1\right]\allowbreak\, \left[2\nv + 3\right]\allowbreak\, \left[2\nv + 2\right]\allowbreak\, \left[\nv + 2\right]\allowbreak\, \left[\nv + 1\right]\allowbreak\, \left[\nv + 1\right]\allowbreak\, \left[\nv\right]\allowbreak\, \left[\nv - 1\right]\allowbreak\, \left[3\right]\allowbreak\, \left[2\right]\allowbreak\, \left[1\right] \\
\DataV{4}{0}{2}{1} & \left[6\nv + 5\right]\allowbreak\, \left[6\nv + 3\right]\allowbreak\, \left[6\nv + 2\right]\allowbreak\, \left[6\nv\right]\allowbreak\, \left[5\nv + 2\right]\allowbreak\, \left[5\nv + 1\right]\allowbreak\, \left[5\nv\right]\allowbreak\, \left[5\nv - 1\right]\allowbreak\, \left[4\nv + 1\right]\allowbreak\, \left[4\nv\right]\allowbreak\, \left[4\nv - 1\right]\allowbreak\, \left[4\nv - 2\right]\allowbreak\, \left[3\nv\right]\allowbreak\, \left[3\nv - 3\right] & \left[3\nv + 4\right]\allowbreak\, \left[3\nv + 1\right]\allowbreak\, \left[2\nv + 3\right]\allowbreak\, \left[2\nv + 2\right]\allowbreak\, \left[2\nv + 1\right]\allowbreak\, \left[2\nv\right]\allowbreak\, \left[\nv + 2\right]\allowbreak\, \left[\nv + 1\right]\allowbreak\, \left[\nv\right]\allowbreak\, \left[\nv - 1\right]\allowbreak\, \left[4\right]\allowbreak\, \left[2\right]\allowbreak\, \left[1\right]\allowbreak\, \left[1\right] \\
\DataV{4}{0}{1}{2} & \left[6\nv + 6\right]\allowbreak\, \left[6\nv + 3\right]\allowbreak\, \left[6\nv + 1\right]\allowbreak\, \left[6\nv\right]\allowbreak\, \left[5\nv + 2\right]\allowbreak\, \left[5\nv + 1\right]\allowbreak\, \left[5\nv\right]\allowbreak\, \left[5\nv - 1\right]\allowbreak\, \left[4\nv + 3\right]\allowbreak\, \left[4\nv + 1\right]\allowbreak\, \left[4\nv\right]\allowbreak\, \left[4\nv - 2\right]\allowbreak\, \left[3\nv - 1\right] & \left[3\nv + 3\right]\allowbreak\, \left[2\nv + 4\right]\allowbreak\, \left[2\nv + 2\right]\allowbreak\, \left[2\nv + 1\right]\allowbreak\, \left[2\nv - 1\right]\allowbreak\, \left[\nv + 3\right]\allowbreak\, \left[\nv + 2\right]\allowbreak\, \left[\nv + 1\right]\allowbreak\, \left[\nv\right]\allowbreak\, \left[4\right]\allowbreak\, \left[2\right]\allowbreak\, \left[1\right]\allowbreak\, \left[1\right] \\
\DataV{4}{1}{3}{2} & \left[6\nv + 4\right]\allowbreak\, \left[6\nv + 3\right]\allowbreak\, \left[6\nv + 1\right]\allowbreak\, \left[6\nv\right]\allowbreak\, \left[5\nv + 2\right]\allowbreak\, \left[5\nv\right]\allowbreak\, \left[5\nv - 1\right]\allowbreak\, \left[4\nv + 1\right]\allowbreak\, \left[4\nv\right]\allowbreak\, \left[4\nv - 2\right]\allowbreak\, \left[3\nv\right]\allowbreak\, \left[3\nv - 1\right]\allowbreak\, \left[3\nv - 3\right] & \left[3\nv + 2\right]\allowbreak\, \left[2\nv + 4\right]\allowbreak\, \left[2\nv + 2\right]\allowbreak\, \left[2\nv + 1\right]\allowbreak\, \left[2\nv\right]\allowbreak\, \left[\nv + 3\right]\allowbreak\, \left[\nv + 2\right]\allowbreak\, \left[\nv + 1\right]\allowbreak\, \left[\nv\right]\allowbreak\, \left[\nv - 1\right]\allowbreak\, \left[2\right]\allowbreak\, \left[1\right]\allowbreak\, \left[1\right] \\
\DataV{4}{1}{2}{3} & \left[6\nv + 5\right]\allowbreak\, \left[6\nv + 2\right]\allowbreak\, \left[6\nv + 1\right]\allowbreak\, \left[6\nv\right]\allowbreak\, \left[5\nv + 2\right]\allowbreak\, \left[5\nv\right]\allowbreak\, \left[5\nv - 1\right]\allowbreak\, \left[4\nv + 2\right]\allowbreak\, \left[4\nv\right]\allowbreak\, \left[4\nv - 1\right]\allowbreak\, \left[3\nv - 1\right]\allowbreak\, \left[3\nv - 3\right] & \left[2\nv + 4\right]\allowbreak\, \left[2\nv + 2\right]\allowbreak\, \left[2\nv + 1\right]\allowbreak\, \left[2\nv\right]\allowbreak\, \left[\nv + 4\right]\allowbreak\, \left[\nv + 2\right]\allowbreak\, \left[\nv + 1\right]\allowbreak\, \left[\nv + 1\right]\allowbreak\, \left[\nv - 1\right]\allowbreak\, \left[2\right]\allowbreak\, \left[1\right]\allowbreak\, \left[1\right] \\
\DataV{4}{3}{2}{1} & -\left[\nv - 5\right]\allowbreak\, \left[4\nv - 2\right]\allowbreak\, \left[5\nv - 1\right]\allowbreak\, \left[6\nv\right]\allowbreak\, \left[3\nv - 2\right]\allowbreak\, \left[5\nv\right]\allowbreak\, \left[6\nv + 1\right]\allowbreak\, \left[2\nv - 2\right]\allowbreak\, \left[4\nv\right]\allowbreak\, \left[5\nv + 1\right]\allowbreak\, \left[4\nv + 1\right]\allowbreak\, \left[6\nv + 3\right] & \left[2\nv + 4\right]\allowbreak\, \left[2\right]\allowbreak\, \left[\nv - 1\right]\allowbreak\, \left[2\nv\right]\allowbreak\, \left[3\nv + 4\right]\allowbreak\, \left[\nv + 2\right]\allowbreak\, \left[1\right]\allowbreak\, \left[1\right]\allowbreak\, \left[2\nv + 1\right]\allowbreak\, \left[2\nv + 2\right]\allowbreak\, \left[\nv + 1\right]\allowbreak\, \left[\nv + 1\right] \\
\DataV{4}{2}{3}{1} & \left[2\nv - 4\right]\allowbreak\, \left[3\nv - 3\right]\allowbreak\, \left[5\nv - 1\right]\allowbreak\, \left[6\nv\right]\allowbreak\, \left[3\nv - 2\right]\allowbreak\, \left[5\nv\right]\allowbreak\, \left[6\nv + 1\right]\allowbreak\, \left[3\nv - 1\right]\allowbreak\, \left[4\nv\right]\allowbreak\, \left[6\nv + 2\right]\allowbreak\, \left[3\nv\right]\allowbreak\, \left[4\nv + 1\right]\allowbreak\, \left[6\nv + 3\right] & \left[\nv - 2\right]\allowbreak\, \left[2\nv + 4\right]\allowbreak\, \left[2\right]\allowbreak\, \left[\nv - 1\right]\allowbreak\, \left[2\nv\right]\allowbreak\, \left[2\nv + 3\right]\allowbreak\, \left[\nv + 2\right]\allowbreak\, \left[1\right]\allowbreak\, \left[\nv\right]\allowbreak\, \left[2\nv + 1\right]\allowbreak\, \left[2\nv + 2\right]\allowbreak\, \left[\nv + 1\right]\allowbreak\, \left[\nv + 1\right] \\
\DataV{4}{2}{4}{2} & \left[2\nv - 1\right]\allowbreak\, \left[3\nv - 1\right]\allowbreak\, \left[4\nv - 1\right]\allowbreak\, \left[5\nv - 1\right]\allowbreak\, \left[3\nv + 1\right]\allowbreak\, \left[4\nv + 1\right]\allowbreak\, \left[5\nv + 1\right]\allowbreak\, \left[6\nv + 1\right]\allowbreak\, \left[2\nv - 2\right]\allowbreak\, \left[2\nv - 4\right]\allowbreak\, \left[2\nv\right]\allowbreak\, \left[6\nv\right]\allowbreak\, \left[4\nv\right]\allowbreak\, \left[4\nv - 2\right]\allowbreak\, \left[4\nv + 2\right]\allowbreak\, \left[6\nv + 2\right]\allowbreak\, \left[6\nv + 4\right] & \left[1\right]\allowbreak\, \left[1\right]\allowbreak\, \left[\nv + 1\right]\allowbreak\, \left[\nv + 1\right]\allowbreak\, \left[\nv + 2\right]\allowbreak\, \left[\nv + 2\right]\allowbreak\, \left[\nv + 4\right]\allowbreak\, \left[3\nv + 4\right]\allowbreak\, \left[\nv - 1\right]\allowbreak\, \left[\nv - 2\right]\allowbreak\, \left[\nv\right]\allowbreak\, \left[3\nv\right]\allowbreak\, \left[2\nv\right]\allowbreak\, \left[2\nv - 1\right]\allowbreak\, \left[2\nv + 1\right]\allowbreak\, \left[3\nv + 1\right]\allowbreak\, \left[3\nv + 2\right] \\
\bottomrule
\end{longtable}
\normalsize

%% file: qdim/level5.tex
\scriptsize
\begin{longtable}{@{}>{$}l<{$} >{$\displaystyle}p{0.39\linewidth}<{$} >{$\displaystyle}p{0.39\linewidth}<{$}@{}}
\toprule
V & N(V) & D(V) \\
\midrule\endhead
\DataV{5}{1}{0}{0} & \left[4\nv\right]\allowbreak\, \left[3\nv - 1\right]\allowbreak\, \left[2\nv - 2\right]\allowbreak\, \left[\nv - 3\right]\allowbreak\, \left[-4\right]\allowbreak\, \left[6\nv + 1\right]\allowbreak\, \left[5\nv\right]\allowbreak\, \left[4\nv - 1\right]\allowbreak\, \left[3\nv - 2\right]\allowbreak\, \left[2\nv - 3\right]\allowbreak\, \left[6\nv\right]\allowbreak\, \left[5\nv - 1\right]\allowbreak\, \left[4\nv - 2\right]\allowbreak\, \left[3\nv - 3\right]\allowbreak\, \left[-3\nv - 9\right] & \left[-\nv - 1\right]\allowbreak\, \left[-2\nv - 2\right]\allowbreak\, \left[-3\nv - 3\right]\allowbreak\, \left[-4\nv - 4\right]\allowbreak\, \left[-5\nv - 5\right]\allowbreak\, \left[-1\right]\allowbreak\, \left[-\nv - 2\right]\allowbreak\, \left[-2\nv - 3\right]\allowbreak\, \left[-3\nv - 4\right]\allowbreak\, \left[-4\nv - 5\right]\allowbreak\, \left[2\nv\right]\allowbreak\, \left[\nv - 1\right]\allowbreak\, \left[-2\right]\allowbreak\, \left[-\nv - 3\right]\allowbreak\, \left[-2\nv - 4\right] \\
\DataV{5}{0}{0}{1} & \left[4\nv\right]\allowbreak\, \left[4\nv + 1\right]\allowbreak\, \left[4\nv + 2\right]\allowbreak\, \left[4\nv + 3\right]\allowbreak\, \left[4\nv + 4\right]\allowbreak\, \left[5\nv - 1\right]\allowbreak\, \left[5\nv\right]\allowbreak\, \left[5\nv + 1\right]\allowbreak\, \left[5\nv + 2\right]\allowbreak\, \left[5\nv + 3\right]\allowbreak\, \left[6\nv\right]\allowbreak\, \left[6\nv + 1\right]\allowbreak\, \left[6\nv + 2\right]\allowbreak\, \left[6\nv + 3\right]\allowbreak\, \left[6\nv + 9\right] & \left[1\right]\allowbreak\, \left[2\right]\allowbreak\, \left[3\right]\allowbreak\, \left[4\right]\allowbreak\, \left[5\right]\allowbreak\, \left[\nv + 1\right]\allowbreak\, \left[\nv + 2\right]\allowbreak\, \left[\nv + 3\right]\allowbreak\, \left[\nv + 4\right]\allowbreak\, \left[\nv + 5\right]\allowbreak\, \left[2\nv\right]\allowbreak\, \left[2\nv + 1\right]\allowbreak\, \left[2\nv + 2\right]\allowbreak\, \left[2\nv + 3\right]\allowbreak\, \left[2\nv + 4\right] \\
\DataV{5}{0}{1}{0} & \wedge^5(L) - \wedge^3(L)\,L + 2\,\wedge^2(L)\,L - 2\,L\,L + L & 1 \\
\DataV{5}{0}{1}{3} & \left[3\nv - 1\right]\allowbreak\, \left[4\nv - 2\right]\allowbreak\, \left[4\nv\right]\allowbreak\, \left[4\nv + 1\right]\allowbreak\, \left[4\nv + 2\right]\allowbreak\, \left[4\nv + 4\right]\allowbreak\, \left[5\nv - 1\right]\allowbreak\, \left[5\nv\right]\allowbreak\, \left[5\nv + 1\right]\allowbreak\, \left[5\nv + 2\right]\allowbreak\, \left[5\nv + 3\right]\allowbreak\, \left[6\nv\right]\allowbreak\, \left[6\nv + 1\right]\allowbreak\, \left[6\nv + 2\right]\allowbreak\, \left[6\nv + 4\right]\allowbreak\, \left[6\nv + 8\right] & \left[1\right]\allowbreak\, \left[1\right]\allowbreak\, \left[2\right]\allowbreak\, \left[3\right]\allowbreak\, \left[5\right]\allowbreak\, \left[\nv\right]\allowbreak\, \left[\nv + 1\right]\allowbreak\, \left[\nv + 2\right]\allowbreak\, \left[\nv + 3\right]\allowbreak\, \left[\nv + 4\right]\allowbreak\, \left[2\nv - 1\right]\allowbreak\, \left[2\nv + 1\right]\allowbreak\, \left[2\nv + 2\right]\allowbreak\, \left[2\nv + 3\right]\allowbreak\, \left[2\nv + 5\right]\allowbreak\, \left[3\nv + 4\right] \\
\DataV{5}{3}{1}{0} & \left[4\nv + 1\right]\allowbreak\, \left[6\nv + 2\right]\allowbreak\, \left[4\nv\right]\allowbreak\, \left[3\nv - 1\right]\allowbreak\, \left[2\nv - 2\right]\allowbreak\, \left[-4\right]\allowbreak\, \left[6\nv + 1\right]\allowbreak\, \left[5\nv\right]\allowbreak\, \left[4\nv - 1\right]\allowbreak\, \left[3\nv - 2\right]\allowbreak\, \left[2\nv - 3\right]\allowbreak\, \left[6\nv\right]\allowbreak\, \left[5\nv - 1\right]\allowbreak\, \left[4\nv - 2\right]\allowbreak\, \left[2\nv - 4\right]\allowbreak\, \left[-2\nv - 8\right] & \left[-\nv - 1\right]\allowbreak\, \left[-\nv - 1\right]\allowbreak\, \left[-2\nv - 2\right]\allowbreak\, \left[-3\nv - 3\right]\allowbreak\, \left[-5\nv - 5\right]\allowbreak\, \left[\nv\right]\allowbreak\, \left[-1\right]\allowbreak\, \left[-\nv - 2\right]\allowbreak\, \left[-2\nv - 3\right]\allowbreak\, \left[-3\nv - 4\right]\allowbreak\, \left[3\nv + 1\right]\allowbreak\, \left[\nv - 1\right]\allowbreak\, \left[-2\right]\allowbreak\, \left[-\nv - 3\right]\allowbreak\, \left[-3\nv - 5\right]\allowbreak\, \left[-\nv - 4\right] \\
\DataV{5}{0}{2}{3} & \left[3\nv - 1\right]\allowbreak\, \left[3\nv\right]\allowbreak\, \left[4\nv - 2\right]\allowbreak\, \left[4\nv - 1\right]\allowbreak\, \left[4\nv\right]\allowbreak\, \left[4\nv + 1\right]\allowbreak\, \left[4\nv + 4\right]\allowbreak\, \left[5\nv - 1\right]\allowbreak\, \left[5\nv\right]\allowbreak\, \left[5\nv + 1\right]\allowbreak\, \left[5\nv + 2\right]\allowbreak\, \left[5\nv + 3\right]\allowbreak\, \left[6\nv\right]\allowbreak\, \left[6\nv + 1\right]\allowbreak\, \left[6\nv + 2\right]\allowbreak\, \left[6\nv + 5\right]\allowbreak\, \left[6\nv + 7\right] & \left[1\right]\allowbreak\, \left[1\right]\allowbreak\, \left[2\right]\allowbreak\, \left[3\right]\allowbreak\, \left[4\right]\allowbreak\, \left[\nv\right]\allowbreak\, \left[\nv + 1\right]\allowbreak\, \left[\nv + 1\right]\allowbreak\, \left[\nv + 2\right]\allowbreak\, \left[\nv + 3\right]\allowbreak\, \left[2\nv - 1\right]\allowbreak\, \left[2\nv\right]\allowbreak\, \left[2\nv + 2\right]\allowbreak\, \left[2\nv + 4\right]\allowbreak\, \left[2\nv + 5\right]\allowbreak\, \left[3\nv + 3\right]\allowbreak\, \left[3\nv + 4\right] \\
\DataV{5}{0}{3}{2} & \left[3\nv - 3\right]\allowbreak\, \left[3\nv - 1\right]\allowbreak\, \left[4\nv - 2\right]\allowbreak\, \left[4\nv - 1\right]\allowbreak\, \left[4\nv\right]\allowbreak\, \left[4\nv\right]\allowbreak\, \left[4\nv + 1\right]\allowbreak\, \left[5\nv - 1\right]\allowbreak\, \left[5\nv\right]\allowbreak\, \left[5\nv + 1\right]\allowbreak\, \left[5\nv + 2\right]\allowbreak\, \left[5\nv + 3\right]\allowbreak\, \left[6\nv\right]\allowbreak\, \left[6\nv + 1\right]\allowbreak\, \left[6\nv + 3\right]\allowbreak\, \left[6\nv + 5\right]\allowbreak\, \left[6\nv + 6\right] & \left[1\right]\allowbreak\, \left[1\right]\allowbreak\, \left[2\right]\allowbreak\, \left[3\right]\allowbreak\, \left[4\right]\allowbreak\, \left[\nv - 1\right]\allowbreak\, \left[\nv\right]\allowbreak\, \left[\nv + 1\right]\allowbreak\, \left[\nv + 1\right]\allowbreak\, \left[\nv + 2\right]\allowbreak\, \left[2\nv\right]\allowbreak\, \left[2\nv + 1\right]\allowbreak\, \left[2\nv + 2\right]\allowbreak\, \left[2\nv + 3\right]\allowbreak\, \left[2\nv + 4\right]\allowbreak\, \left[3\nv + 3\right]\allowbreak\, \left[3\nv + 5\right] \\
\DataV{5}{2}{3}{0} & \left[6\nv + 3\right]\allowbreak\, \left[4\nv + 1\right]\allowbreak\, \left[6\nv + 2\right]\allowbreak\, \left[5\nv + 1\right]\allowbreak\, \left[4\nv\right]\allowbreak\, \left[4\nv\right]\allowbreak\, \left[3\nv - 1\right]\allowbreak\, \left[6\nv + 1\right]\allowbreak\, \left[5\nv\right]\allowbreak\, \left[4\nv - 1\right]\allowbreak\, \left[3\nv - 2\right]\allowbreak\, \left[2\nv - 3\right]\allowbreak\, \left[6\nv\right]\allowbreak\, \left[5\nv - 1\right]\allowbreak\, \left[3\nv - 3\right]\allowbreak\, \left[\nv - 5\right]\allowbreak\, \left[-6\right] & \left[-\nv - 1\right]\allowbreak\, \left[-\nv - 1\right]\allowbreak\, \left[-2\nv - 2\right]\allowbreak\, \left[-3\nv - 3\right]\allowbreak\, \left[-4\nv - 4\right]\allowbreak\, \left[2\nv + 1\right]\allowbreak\, \left[\nv\right]\allowbreak\, \left[-1\right]\allowbreak\, \left[-1\right]\allowbreak\, \left[-\nv - 2\right]\allowbreak\, \left[2\nv\right]\allowbreak\, \left[\nv - 1\right]\allowbreak\, \left[-2\right]\allowbreak\, \left[-\nv - 3\right]\allowbreak\, \left[-2\nv - 4\right]\allowbreak\, \left[-3\right]\allowbreak\, \left[-2\nv - 5\right] \\
\DataV{5}{3}{2}{0} & \left[4\nv + 1\right]\allowbreak\, \left[3\nv\right]\allowbreak\, \left[6\nv + 2\right]\allowbreak\, \left[5\nv + 1\right]\allowbreak\, \left[4\nv\right]\allowbreak\, \left[3\nv - 1\right]\allowbreak\, \left[-4\right]\allowbreak\, \left[6\nv + 1\right]\allowbreak\, \left[5\nv\right]\allowbreak\, \left[4\nv - 1\right]\allowbreak\, \left[3\nv - 2\right]\allowbreak\, \left[2\nv - 3\right]\allowbreak\, \left[6\nv\right]\allowbreak\, \left[5\nv - 1\right]\allowbreak\, \left[4\nv - 2\right]\allowbreak\, \left[\nv - 5\right]\allowbreak\, \left[-\nv - 7\right] & \left[-\nv - 1\right]\allowbreak\, \left[-\nv - 1\right]\allowbreak\, \left[-2\nv - 2\right]\allowbreak\, \left[-3\nv - 3\right]\allowbreak\, \left[-4\nv - 4\right]\allowbreak\, \left[\nv\right]\allowbreak\, \left[-1\right]\allowbreak\, \left[-1\right]\allowbreak\, \left[-\nv - 2\right]\allowbreak\, \left[-2\nv - 3\right]\allowbreak\, \left[3\nv + 1\right]\allowbreak\, \left[2\nv\right]\allowbreak\, \left[-2\right]\allowbreak\, \left[-2\nv - 4\right]\allowbreak\, \left[-3\nv - 5\right]\allowbreak\, \left[-3\right]\allowbreak\, \left[-\nv - 4\right] \\
\DataV{5}{0}{3}{3} & \left[3\nv - 3\right]\allowbreak\, \left[3\nv\right]\allowbreak\, \left[4\nv - 2\right]\allowbreak\, \left[4\nv - 1\right]\allowbreak\, \left[4\nv\right]\allowbreak\, \left[4\nv + 1\right]\allowbreak\, \left[4\nv + 2\right]\allowbreak\, \left[5\nv - 1\right]\allowbreak\, \left[5\nv\right]\allowbreak\, \left[5\nv + 1\right]\allowbreak\, \left[5\nv + 2\right]\allowbreak\, \left[5\nv + 3\right]\allowbreak\, \left[6\nv\right]\allowbreak\, \left[6\nv + 1\right]\allowbreak\, \left[6\nv + 3\right]\allowbreak\, \left[6\nv + 4\right]\allowbreak\, \left[6\nv + 7\right] & \left[1\right]\allowbreak\, \left[1\right]\allowbreak\, \left[2\right]\allowbreak\, \left[2\right]\allowbreak\, \left[5\right]\allowbreak\, \left[\nv - 1\right]\allowbreak\, \left[\nv\right]\allowbreak\, \left[\nv + 1\right]\allowbreak\, \left[\nv + 2\right]\allowbreak\, \left[\nv + 3\right]\allowbreak\, \left[2\nv\right]\allowbreak\, \left[2\nv + 1\right]\allowbreak\, \left[2\nv + 2\right]\allowbreak\, \left[2\nv + 3\right]\allowbreak\, \left[2\nv + 4\right]\allowbreak\, \left[3\nv + 2\right]\allowbreak\, \left[3\nv + 5\right] \\
\DataV{5}{3}{3}{0} & \left[6\nv + 3\right]\allowbreak\, \left[3\nv\right]\allowbreak\, \left[6\nv + 2\right]\allowbreak\, \left[5\nv + 1\right]\allowbreak\, \left[4\nv\right]\allowbreak\, \left[3\nv - 1\right]\allowbreak\, \left[2\nv - 2\right]\allowbreak\, \left[6\nv + 1\right]\allowbreak\, \left[5\nv\right]\allowbreak\, \left[4\nv - 1\right]\allowbreak\, \left[3\nv - 2\right]\allowbreak\, \left[2\nv - 3\right]\allowbreak\, \left[6\nv\right]\allowbreak\, \left[5\nv - 1\right]\allowbreak\, \left[3\nv - 3\right]\allowbreak\, \left[2\nv - 4\right]\allowbreak\, \left[-\nv - 7\right] & \left[-\nv - 1\right]\allowbreak\, \left[-\nv - 1\right]\allowbreak\, \left[-2\nv - 2\right]\allowbreak\, \left[-2\nv - 2\right]\allowbreak\, \left[-5\nv - 5\right]\allowbreak\, \left[2\nv + 1\right]\allowbreak\, \left[\nv\right]\allowbreak\, \left[-1\right]\allowbreak\, \left[-\nv - 2\right]\allowbreak\, \left[-2\nv - 3\right]\allowbreak\, \left[2\nv\right]\allowbreak\, \left[\nv - 1\right]\allowbreak\, \left[-2\right]\allowbreak\, \left[-\nv - 3\right]\allowbreak\, \left[-2\nv - 4\right]\allowbreak\, \left[\nv - 2\right]\allowbreak\, \left[-2\nv - 5\right] \\
\DataV{5}{1}{6}{3} & \left[2\nv - 2\right]\allowbreak\, \left[3\nv - 3\right]\allowbreak\, \left[3\nv - 2\right]\allowbreak\, \left[3\nv - 1\right]\allowbreak\, \left[4\nv - 2\right]\allowbreak\, \left[4\nv - 1\right]\allowbreak\, \left[4\nv\right]\allowbreak\, \left[4\nv + 1\right]\allowbreak\, \left[4\nv + 3\right]\allowbreak\, \left[5\nv - 1\right]\allowbreak\, \left[5\nv\right]\allowbreak\, \left[5\nv + 1\right]\allowbreak\, \left[5\nv + 3\right]\allowbreak\, \left[6\nv\right]\allowbreak\, \left[6\nv + 1\right]\allowbreak\, \left[6\nv + 3\right]\allowbreak\, \left[6\nv + 4\right]\allowbreak\, \left[6\nv + 5\right] & \left[1\right]\allowbreak\, \left[1\right]\allowbreak\, \left[2\right]\allowbreak\, \left[3\right]\allowbreak\, \left[\nv - 1\right]\allowbreak\, \left[\nv - 1\right]\allowbreak\, \left[\nv + 1\right]\allowbreak\, \left[\nv + 1\right]\allowbreak\, \left[\nv + 2\right]\allowbreak\, \left[\nv + 3\right]\allowbreak\, \left[2\nv - 1\right]\allowbreak\, \left[2\nv\right]\allowbreak\, \left[2\nv + 1\right]\allowbreak\, \left[2\nv + 2\right]\allowbreak\, \left[2\nv + 3\right]\allowbreak\, \left[2\nv + 4\right]\allowbreak\, \left[3\nv + 4\right]\allowbreak\, \left[4\nv + 4\right] \\
\DataV{5}{3}{6}{1} & \left[4\nv + 2\right]\allowbreak\, \left[6\nv + 3\right]\allowbreak\, \left[5\nv + 2\right]\allowbreak\, \left[4\nv + 1\right]\allowbreak\, \left[6\nv + 2\right]\allowbreak\, \left[5\nv + 1\right]\allowbreak\, \left[4\nv\right]\allowbreak\, \left[3\nv - 1\right]\allowbreak\, \left[\nv - 3\right]\allowbreak\, \left[6\nv + 1\right]\allowbreak\, \left[5\nv\right]\allowbreak\, \left[4\nv - 1\right]\allowbreak\, \left[2\nv - 3\right]\allowbreak\, \left[6\nv\right]\allowbreak\, \left[5\nv - 1\right]\allowbreak\, \left[3\nv - 3\right]\allowbreak\, \left[2\nv - 4\right]\allowbreak\, \left[\nv - 5\right] & \left[-\nv - 1\right]\allowbreak\, \left[-\nv - 1\right]\allowbreak\, \left[-2\nv - 2\right]\allowbreak\, \left[-3\nv - 3\right]\allowbreak\, \left[2\nv + 1\right]\allowbreak\, \left[2\nv + 1\right]\allowbreak\, \left[-1\right]\allowbreak\, \left[-1\right]\allowbreak\, \left[-\nv - 2\right]\allowbreak\, \left[-2\nv - 3\right]\allowbreak\, \left[3\nv + 1\right]\allowbreak\, \left[2\nv\right]\allowbreak\, \left[\nv - 1\right]\allowbreak\, \left[-2\right]\allowbreak\, \left[-\nv - 3\right]\allowbreak\, \left[-2\nv - 4\right]\allowbreak\, \left[-\nv - 4\right]\allowbreak\, \left[-4\right] \\
\DataV{5}{1}{3}{6} & \left[3\nv - 3\right]\allowbreak\, \left[3\nv - 1\right]\allowbreak\, \left[4\nv - 1\right]\allowbreak\, \left[4\nv\right]\allowbreak\, \left[4\nv + 1\right]\allowbreak\, \left[4\nv + 3\right]\allowbreak\, \left[5\nv - 1\right]\allowbreak\, \left[5\nv\right]\allowbreak\, \left[5\nv + 1\right]\allowbreak\, \left[5\nv + 3\right]\allowbreak\, \left[6\nv\right]\allowbreak\, \left[6\nv + 1\right]\allowbreak\, \left[6\nv + 2\right]\allowbreak\, \left[6\nv + 3\right]\allowbreak\, \left[6\nv + 7\right] & \left[1\right]\allowbreak\, \left[1\right]\allowbreak\, \left[2\right]\allowbreak\, \left[3\right]\allowbreak\, \left[\nv - 1\right]\allowbreak\, \left[\nv + 1\right]\allowbreak\, \left[\nv + 1\right]\allowbreak\, \left[\nv + 2\right]\allowbreak\, \left[\nv + 3\right]\allowbreak\, \left[\nv + 5\right]\allowbreak\, \left[2\nv\right]\allowbreak\, \left[2\nv + 1\right]\allowbreak\, \left[2\nv + 2\right]\allowbreak\, \left[2\nv + 3\right]\allowbreak\, \left[2\nv + 5\right] \\
\DataV{5}{6}{3}{1} & \left[6\nv + 3\right]\allowbreak\, \left[4\nv + 1\right]\allowbreak\, \left[5\nv + 1\right]\allowbreak\, \left[4\nv\right]\allowbreak\, \left[3\nv - 1\right]\allowbreak\, \left[\nv - 3\right]\allowbreak\, \left[6\nv + 1\right]\allowbreak\, \left[5\nv\right]\allowbreak\, \left[4\nv - 1\right]\allowbreak\, \left[2\nv - 3\right]\allowbreak\, \left[6\nv\right]\allowbreak\, \left[5\nv - 1\right]\allowbreak\, \left[4\nv - 2\right]\allowbreak\, \left[3\nv - 3\right]\allowbreak\, \left[-\nv - 7\right] & \left[-\nv - 1\right]\allowbreak\, \left[-\nv - 1\right]\allowbreak\, \left[-2\nv - 2\right]\allowbreak\, \left[-3\nv - 3\right]\allowbreak\, \left[2\nv + 1\right]\allowbreak\, \left[-1\right]\allowbreak\, \left[-1\right]\allowbreak\, \left[-\nv - 2\right]\allowbreak\, \left[-2\nv - 3\right]\allowbreak\, \left[-4\nv - 5\right]\allowbreak\, \left[2\nv\right]\allowbreak\, \left[\nv - 1\right]\allowbreak\, \left[-2\right]\allowbreak\, \left[-\nv - 3\right]\allowbreak\, \left[-3\nv - 5\right] \\
\DataV{5}{3}{6}{6} & \left[3\nv - 3\right]\allowbreak\, \left[3\nv - 2\right]\allowbreak\, \left[3\nv - 1\right]\allowbreak\, \left[3\nv\right]\allowbreak\, \left[4\nv - 2\right]\allowbreak\, \left[4\nv\right]\allowbreak\, \left[4\nv + 1\right]\allowbreak\, \left[5\nv - 1\right]\allowbreak\, \left[5\nv\right]\allowbreak\, \left[5\nv + 3\right]\allowbreak\, \left[6\nv\right]\allowbreak\, \left[6\nv + 1\right]\allowbreak\, \left[6\nv + 2\right]\allowbreak\, \left[6\nv + 3\right]\allowbreak\, \left[6\nv + 5\right] & \left[1\right]\allowbreak\, \left[1\right]\allowbreak\, \left[2\right]\allowbreak\, \left[\nv - 1\right]\allowbreak\, \left[\nv\right]\allowbreak\, \left[\nv + 1\right]\allowbreak\, \left[\nv + 1\right]\allowbreak\, \left[\nv + 2\right]\allowbreak\, \left[\nv + 3\right]\allowbreak\, \left[\nv + 4\right]\allowbreak\, \left[2\nv\right]\allowbreak\, \left[2\nv + 1\right]\allowbreak\, \left[2\nv + 2\right]\allowbreak\, \left[2\nv + 3\right]\allowbreak\, \left[2\nv + 4\right] \\
\DataV{5}{6}{6}{3} & \left[6\nv + 3\right]\allowbreak\, \left[5\nv + 2\right]\allowbreak\, \left[4\nv + 1\right]\allowbreak\, \left[3\nv\right]\allowbreak\, \left[6\nv + 2\right]\allowbreak\, \left[4\nv\right]\allowbreak\, \left[3\nv - 1\right]\allowbreak\, \left[6\nv + 1\right]\allowbreak\, \left[5\nv\right]\allowbreak\, \left[2\nv - 3\right]\allowbreak\, \left[6\nv\right]\allowbreak\, \left[5\nv - 1\right]\allowbreak\, \left[4\nv - 2\right]\allowbreak\, \left[3\nv - 3\right]\allowbreak\, \left[\nv - 5\right] & \left[-\nv - 1\right]\allowbreak\, \left[-\nv - 1\right]\allowbreak\, \left[-2\nv - 2\right]\allowbreak\, \left[2\nv + 1\right]\allowbreak\, \left[\nv\right]\allowbreak\, \left[-1\right]\allowbreak\, \left[-1\right]\allowbreak\, \left[-\nv - 2\right]\allowbreak\, \left[-2\nv - 3\right]\allowbreak\, \left[-3\nv - 4\right]\allowbreak\, \left[2\nv\right]\allowbreak\, \left[\nv - 1\right]\allowbreak\, \left[-2\right]\allowbreak\, \left[-\nv - 3\right]\allowbreak\, \left[-2\nv - 4\right] \\
\DataV{5}{2}{9}{9} & \left[3\nv - 3\right]\allowbreak\, \left[3\nv - 1\right]\allowbreak\, \left[3\nv\right]\allowbreak\, \left[4\nv - 2\right]\allowbreak\, \left[4\nv\right]\allowbreak\, \left[4\nv\right]\allowbreak\, \left[4\nv + 2\right]\allowbreak\, \left[5\nv - 1\right]\allowbreak\, \left[5\nv\right]\allowbreak\, \left[5\nv + 1\right]\allowbreak\, \left[5\nv + 3\right]\allowbreak\, \left[6\nv\right]\allowbreak\, \left[6\nv + 1\right]\allowbreak\, \left[6\nv + 2\right]\allowbreak\, \left[6\nv + 4\right]\allowbreak\, \left[6\nv + 6\right] & \left[1\right]\allowbreak\, \left[1\right]\allowbreak\, \left[1\right]\allowbreak\, \left[3\right]\allowbreak\, \left[\nv - 1\right]\allowbreak\, \left[\nv\right]\allowbreak\, \left[\nv + 1\right]\allowbreak\, \left[\nv + 2\right]\allowbreak\, \left[\nv + 2\right]\allowbreak\, \left[\nv + 4\right]\allowbreak\, \left[2\nv\right]\allowbreak\, \left[2\nv + 1\right]\allowbreak\, \left[2\nv + 2\right]\allowbreak\, \left[2\nv + 3\right]\allowbreak\, \left[2\nv + 5\right]\allowbreak\, \left[3\nv + 3\right] \\
\DataV{5}{9}{9}{2} & \left[6\nv + 3\right]\allowbreak\, \left[4\nv + 1\right]\allowbreak\, \left[3\nv\right]\allowbreak\, \left[6\nv + 2\right]\allowbreak\, \left[4\nv\right]\allowbreak\, \left[4\nv\right]\allowbreak\, \left[2\nv - 2\right]\allowbreak\, \left[6\nv + 1\right]\allowbreak\, \left[5\nv\right]\allowbreak\, \left[4\nv - 1\right]\allowbreak\, \left[2\nv - 3\right]\allowbreak\, \left[6\nv\right]\allowbreak\, \left[5\nv - 1\right]\allowbreak\, \left[4\nv - 2\right]\allowbreak\, \left[2\nv - 4\right]\allowbreak\, \left[-6\right] & \left[-\nv - 1\right]\allowbreak\, \left[-\nv - 1\right]\allowbreak\, \left[-\nv - 1\right]\allowbreak\, \left[-3\nv - 3\right]\allowbreak\, \left[2\nv + 1\right]\allowbreak\, \left[\nv\right]\allowbreak\, \left[-1\right]\allowbreak\, \left[-\nv - 2\right]\allowbreak\, \left[-\nv - 2\right]\allowbreak\, \left[-3\nv - 4\right]\allowbreak\, \left[2\nv\right]\allowbreak\, \left[\nv - 1\right]\allowbreak\, \left[-2\right]\allowbreak\, \left[-\nv - 3\right]\allowbreak\, \left[-3\nv - 5\right]\allowbreak\, \left[-3\right] \\
\bottomrule
\end{longtable}
\normalsize
We now discuss the quantum dimension of $V_5(0,1,0)$. From \cite[(8)]{Westbury2006a},
the representation $V_5(0,1,0)$ is
\[ \largewedge^5(L) - \largewedge^3(L)\,L + 2\,\largewedge^2(L)\,L - 2\,L\,L + L \]
Converting from elementary symmetric functions to power sum symmetric functions
\begin{multline*}
	e_{1} - 2 e_{1,1} + 2 e_{2,1} - e_{3,1} + e_{5} = 
p_{1} - 2 p_{1,1} + p_{1,1,1} - \frac{1}{6} p_{1,1,1,1} + \frac{1}{120} p_{1,1,1,1,1} - p_{2,1} \\+ \frac{1}{2} p_{2,1,1} - \frac{1}{12} p_{2,1,1,1} + \frac{1}{8} p_{2,2,1} - \frac{1}{3} p_{3,1} + \frac{1}{6} p_{3,1,1} - \frac{1}{6} p_{3,2} - \frac{1}{4} p_{4,1} + \frac{1}{5} p_{5}
\end{multline*}
For $k\geqslant 0$, $p_k$ acts by the Adams operator $\psi_k$. Now use
\[ \psi_k(\dim_q(L)) = \frac{[6k\nv+k][5k\nv-k][4k\nv]}{[2k\nv][k\nv+k][k]}    \]
This gives $\dim_q(V_5(0,1,0)) \in \bQ(Q,q)$. It is clear that the denominator is a product
of quantum integers. However the numerator has a large irreducible factor.

The exact rational functions for $\DataV{5}{1}{3}{0}$, $\DataV{5}{0}{3}{1}$, $\DataV{5}{3}{9}{3}$, $\DataV{5}{6}{8}{6}$, $\DataV{5}{3}{9}{8}$, $\DataV{5}{8}{9}{3}$, $\DataV{5}{3}{11}{6}$, $\DataV{5}{6}{11}{3}$, $\DataV{5}{6}{12}{6}$, $\DataV{5}{6}{15}{9}$, $\DataV{5}{9}{15}{6}$ are present in the archive but omitted here.

%% file: orthogonal/degree1.tex
\subsubsection*{Orthogonal functors of degree 1}
\begin{align*}
o_{(1)}(L) &= L.
\end{align*}

%% file: orthogonal/degree2.tex
\subsubsection*{Orthogonal functors of degree 2}
\begin{align*}
o_{(2)}(L) &= \DataV{2}{0}{0}{1} + \DataV{2}{1}{0}{0}.
\end{align*}
\begin{align*}
o_{(1,1)}(L) &= L + \DataV{2}{0}{1}{0}.
\end{align*}

%% file: orthogonal/degree3.tex
\subsubsection*{Orthogonal functors of degree 3}
\begin{align*}
o_{(3)}(L) &= \DataV{2}{0}{1}{0} + \DataV{3}{0}{0}{1} + \DataV{3}{1}{0}{0} + \DataV{3}{1}{1}{1}.
\end{align*}
\begin{align*}
o_{(2,1)}(L) &= L + \DataV{2}{0}{0}{1} + \DataV{2}{0}{1}{0} + \DataV{2}{1}{0}{0} \\
&\quad {}+ \DataV{3}{0}{1}{1} + \DataV{3}{1}{1}{0} + \DataV{3}{1}{1}{1}.
\end{align*}
\begin{align*}
o_{(1,1,1)}(L) &= \DataUnit + \DataV{2}{0}{0}{1} + \DataV{2}{0}{1}{0} + \DataV{2}{1}{0}{0} \\
&\quad {}+ \DataV{3}{0}{1}{0}.
\end{align*}

%% file: orthogonal/degree4.tex
\subsubsection*{Orthogonal functors of degree 4}
\begin{align*}
o_{(4)}(L) &= \DataV{2}{0}{0}{1} + \DataV{2}{1}{0}{0} + \DataV{3}{0}{1}{0} + \DataV{3}{0}{1}{1} \\
&\quad {}+ \DataV{3}{1}{1}{0} + \DataV{4}{0}{0}{1} + \DataV{4}{1}{0}{0} + \DataV{4}{1}{1}{1} \\
&\quad {}+ \DataV{4}{1}{2}{3} + \DataV{4}{3}{2}{1}.
\end{align*}
\begin{align*}
o_{(3,1)}(L) &= L + \DataV{2}{0}{0}{1} + 2\DataV{2}{0}{1}{0} + \DataV{2}{1}{0}{0} \\
&\quad {}+ \DataV{3}{0}{0}{1} + \DataV{3}{0}{1}{0} + 2\DataV{3}{0}{1}{1} + \DataV{3}{1}{0}{0} \\
&\quad {}+ 2\DataV{3}{1}{1}{0} + 3\DataV{3}{1}{1}{1} + \DataV{4}{0}{1}{2} + \DataV{4}{1}{2}{3} \\
&\quad {}+ \DataV{4}{1}{3}{2} + \DataV{4}{2}{1}{0} + \DataV{4}{2}{3}{1} + \DataV{4}{2}{4}{2} \\
&\quad {}+ \DataV{4}{3}{2}{1}.
\end{align*}
\begin{align*}
o_{(2,2)}(L) &= \DataUnit + 2\DataV{2}{0}{0}{1} + \DataV{2}{0}{1}{0} + 2\DataV{2}{1}{0}{0} \\
&\quad {}+ 2\DataV{3}{0}{1}{0} + \DataV{3}{0}{1}{1} + \DataV{3}{1}{1}{0} + \DataV{3}{1}{1}{1} \\
&\quad {}+ \DataV{4}{0}{1}{1} + \DataV{4}{1}{1}{0} + \DataV{4}{1}{1}{1} + \DataV{4}{1}{2}{3} \\
&\quad {}+ \DataV{4}{2}{4}{2} + \DataV{4}{3}{2}{1}.
\end{align*}
\begin{align*}
o_{(2,1,1)}(L) &= 2L + \DataV{2}{0}{0}{1} + 3\DataV{2}{0}{1}{0} + \DataV{2}{1}{0}{0} \\
&\quad {}+ \DataV{3}{0}{0}{1} + \DataV{3}{0}{1}{0} + 2\DataV{3}{0}{1}{1} + \DataV{3}{1}{0}{0} \\
&\quad {}+ 2\DataV{3}{1}{1}{0} + 3\DataV{3}{1}{1}{1} + \DataV{4}{0}{2}{1} + \DataV{4}{1}{2}{0} \\
&\quad {}+ \DataV{4}{1}{3}{2} + \DataV{4}{2}{3}{1} + \DataV{4}{2}{4}{2}.
\end{align*}
\begin{align*}
o_{(1,1,1,1)}(L) &= L + \DataV{2}{0}{0}{1} + \DataV{2}{1}{0}{0} + \DataV{3}{0}{1}{0} \\
&\quad {}+ \DataV{3}{0}{1}{1} + \DataV{3}{1}{1}{0} + \DataV{3}{1}{1}{1} + \DataV{4}{0}{1}{0}.
\end{align*}

%% file: orthogonal/degree5.tex
\subsubsection*{Orthogonal functors of degree 5}
\begin{align*}
o_{(5)}(L) &= L + \DataV{2}{0}{1}{0} + \DataV{3}{0}{0}{1} + \DataV{3}{0}{1}{1} \\
&\quad {}+ \DataV{3}{1}{0}{0} + \DataV{3}{1}{1}{0} + 2\DataV{3}{1}{1}{1} + \DataV{4}{0}{1}{2} \\
&\quad {}+ \DataV{4}{0}{2}{1} + \DataV{4}{1}{2}{0} + \DataV{4}{1}{3}{2} + \DataV{4}{2}{1}{0} \\
&\quad {}+ \DataV{4}{2}{3}{1} + \DataV{4}{2}{4}{2} + \DataV{5}{0}{0}{1} + \DataV{5}{1}{0}{0} \\
&\quad {}+ \DataV{5}{1}{3}{6} + \DataV{5}{3}{6}{6} + \DataV{5}{6}{3}{1} + \DataV{5}{6}{6}{3} \\
&\quad {}+ \DataV{5}{6}{8}{6}.
\end{align*}
\begin{align*}
o_{(4,1)}(L) &= L + 2\DataV{2}{0}{0}{1} + 3\DataV{2}{0}{1}{0} + 2\DataV{2}{1}{0}{0} \\
&\quad {}+ 2\DataV{3}{0}{0}{1} + 3\DataV{3}{0}{1}{0} + 4\DataV{3}{0}{1}{1} + 2\DataV{3}{1}{0}{0} \\
&\quad {}+ 4\DataV{3}{1}{1}{0} + 5\DataV{3}{1}{1}{1} + \DataV{4}{0}{0}{1} + \DataV{4}{0}{1}{0} \\
&\quad {}+ \DataV{4}{0}{1}{1} + 2\DataV{4}{0}{1}{2} + 2\DataV{4}{0}{2}{1} + \DataV{4}{1}{0}{0} \\
&\quad {}+ \DataV{4}{1}{1}{0} + \DataV{4}{1}{1}{1} + 2\DataV{4}{1}{2}{0} + 3\DataV{4}{1}{2}{3} \\
&\quad {}+ 3\DataV{4}{1}{3}{2} + 2\DataV{4}{2}{1}{0} + 3\DataV{4}{2}{3}{1} + 4\DataV{4}{2}{4}{2} \\
&\quad {}+ 3\DataV{4}{3}{2}{1} + \DataV{5}{0}{1}{3} + \DataV{5}{1}{3}{6} + \DataV{5}{2}{9}{9} \\
&\quad {}+ \DataV{5}{3}{1}{0} + \DataV{5}{3}{6}{6} + \DataV{5}{3}{9}{8} + \DataV{5}{6}{12}{6} \\
&\quad {}+ \DataV{5}{6}{15}{9} + \DataV{5}{6}{3}{1} + \DataV{5}{6}{6}{3} + \DataV{5}{6}{8}{6} \\
&\quad {}+ \DataV{5}{8}{9}{3} + \DataV{5}{9}{15}{6} + \DataV{5}{9}{9}{2}.
\end{align*}
\begin{align*}
o_{(3,2)}(L) &= 2L + 2\DataV{2}{0}{0}{1} + 4\DataV{2}{0}{1}{0} + 2\DataV{2}{1}{0}{0} \\
&\quad {}+ 2\DataV{3}{0}{0}{1} + 3\DataV{3}{0}{1}{0} + 5\DataV{3}{0}{1}{1} + 2\DataV{3}{1}{0}{0} \\
&\quad {}+ 5\DataV{3}{1}{1}{0} + 7\DataV{3}{1}{1}{1} + \DataV{4}{0}{1}{0} + \DataV{4}{0}{1}{1} \\
&\quad {}+ 2\DataV{4}{0}{1}{2} + 3\DataV{4}{0}{2}{1} + \DataV{4}{1}{1}{0} + \DataV{4}{1}{1}{1} \\
&\quad {}+ 3\DataV{4}{1}{2}{0} + 3\DataV{4}{1}{2}{3} + 4\DataV{4}{1}{3}{2} + 2\DataV{4}{2}{1}{0} \\
&\quad {}+ 4\DataV{4}{2}{3}{1} + 5\DataV{4}{2}{4}{2} + 3\DataV{4}{3}{2}{1} + \DataV{5}{0}{2}{3} \\
&\quad {}+ \DataV{5}{1}{3}{6} + \DataV{5}{2}{9}{9} + \DataV{5}{3}{11}{6} + \DataV{5}{3}{2}{0} \\
&\quad {}+ \DataV{5}{3}{6}{6} + \DataV{5}{3}{9}{8} + \DataV{5}{6}{11}{3} + \DataV{5}{6}{12}{6} \\
&\quad {}+ \DataV{5}{6}{15}{9} + \DataV{5}{6}{3}{1} + \DataV{5}{6}{6}{3} + 2\DataV{5}{6}{8}{6} \\
&\quad {}+ \DataV{5}{8}{9}{3} + \DataV{5}{9}{15}{6} + \DataV{5}{9}{9}{2}.
\end{align*}
\begin{align*}
o_{(3,1,1)}(L) &= \DataUnit + L + 5\DataV{2}{0}{0}{1} + 4\DataV{2}{0}{1}{0} \\
&\quad {}+ 5\DataV{2}{1}{0}{0} + \DataV{3}{0}{0}{1} + 7\DataV{3}{0}{1}{0} + 6\DataV{3}{0}{1}{1} \\
&\quad {}+ \DataV{3}{1}{0}{0} + 6\DataV{3}{1}{1}{0} + 6\DataV{3}{1}{1}{1} + \DataV{4}{0}{0}{1} \\
&\quad {}+ 2\DataV{4}{0}{1}{0} + 2\DataV{4}{0}{1}{1} + 2\DataV{4}{0}{1}{2} + 2\DataV{4}{0}{2}{1} \\
&\quad {}+ \DataV{4}{1}{0}{0} + 2\DataV{4}{1}{1}{0} + 3\DataV{4}{1}{1}{1} + 2\DataV{4}{1}{2}{0} \\
&\quad {}+ 5\DataV{4}{1}{2}{3} + 4\DataV{4}{1}{3}{2} + 2\DataV{4}{2}{1}{0} + 4\DataV{4}{2}{3}{1} \\
&\quad {}+ 6\DataV{4}{2}{4}{2} + 5\DataV{4}{3}{2}{1} + \DataV{5}{0}{3}{3} + \DataV{5}{1}{6}{3} \\
&\quad {}+ \DataV{5}{2}{9}{9} + \DataV{5}{3}{11}{6} + \DataV{5}{3}{3}{0} + \DataV{5}{3}{6}{1} \\
&\quad {}+ \DataV{5}{3}{9}{3} + \DataV{5}{3}{9}{8} + \DataV{5}{6}{11}{3} + \DataV{5}{6}{12}{6} \\
&\quad {}+ 2\DataV{5}{6}{15}{9} + \DataV{5}{8}{9}{3} + 2\DataV{5}{9}{15}{6} + \DataV{5}{9}{9}{2}.
\end{align*}
\begin{align*}
o_{(2,2,1)}(L) &= 3L + 2\DataV{2}{0}{0}{1} + 5\DataV{2}{0}{1}{0} + 2\DataV{2}{1}{0}{0} \\
&\quad {}+ 2\DataV{3}{0}{0}{1} + 4\DataV{3}{0}{1}{0} + 5\DataV{3}{0}{1}{1} + 2\DataV{3}{1}{0}{0} \\
&\quad {}+ 5\DataV{3}{1}{1}{0} + 7\DataV{3}{1}{1}{1} + 2\DataV{4}{0}{1}{0} + \DataV{4}{0}{1}{1} \\
&\quad {}+ \DataV{4}{0}{1}{2} + 3\DataV{4}{0}{2}{1} + \DataV{4}{1}{1}{0} + \DataV{4}{1}{1}{1} \\
&\quad {}+ 3\DataV{4}{1}{2}{0} + 2\DataV{4}{1}{2}{3} + 4\DataV{4}{1}{3}{2} + \DataV{4}{2}{1}{0} \\
&\quad {}+ 4\DataV{4}{2}{3}{1} + 5\DataV{4}{2}{4}{2} + 2\DataV{4}{3}{2}{1} + \DataV{5}{0}{3}{2} \\
&\quad {}+ \DataV{5}{2}{3}{0} + \DataV{5}{2}{9}{9} + \DataV{5}{3}{11}{6} + \DataV{5}{3}{6}{6} \\
&\quad {}+ \DataV{5}{3}{9}{3} + \DataV{5}{3}{9}{8} + \DataV{5}{6}{11}{3} + \DataV{5}{6}{12}{6} \\
&\quad {}+ \DataV{5}{6}{15}{9} + \DataV{5}{6}{6}{3} + \DataV{5}{6}{8}{6} + \DataV{5}{8}{9}{3} \\
&\quad {}+ \DataV{5}{9}{15}{6} + \DataV{5}{9}{9}{2}.
\end{align*}
\begin{align*}
o_{(2,1,1,1)}(L) &= 2L + 3\DataV{2}{0}{0}{1} + 4\DataV{2}{0}{1}{0} + 3\DataV{2}{1}{0}{0} \\
&\quad {}+ \DataV{3}{0}{0}{1} + 4\DataV{3}{0}{1}{0} + 4\DataV{3}{0}{1}{1} + \DataV{3}{1}{0}{0} \\
&\quad {}+ 4\DataV{3}{1}{1}{0} + 5\DataV{3}{1}{1}{1} + 2\DataV{4}{0}{1}{0} + \DataV{4}{0}{1}{1} \\
&\quad {}+ \DataV{4}{0}{1}{2} + 2\DataV{4}{0}{2}{1} + \DataV{4}{1}{1}{0} + \DataV{4}{1}{1}{1} \\
&\quad {}+ 2\DataV{4}{1}{2}{0} + 2\DataV{4}{1}{2}{3} + 3\DataV{4}{1}{3}{2} + \DataV{4}{2}{1}{0} \\
&\quad {}+ 3\DataV{4}{2}{3}{1} + 4\DataV{4}{2}{4}{2} + 2\DataV{4}{3}{2}{1} + \DataV{5}{0}{3}{1} \\
&\quad {}+ \DataV{5}{1}{3}{0} + \DataV{5}{1}{6}{3} + \DataV{5}{3}{11}{6} + \DataV{5}{3}{6}{1} \\
&\quad {}+ \DataV{5}{3}{9}{3} + \DataV{5}{6}{11}{3} + \DataV{5}{6}{12}{6} + \DataV{5}{6}{15}{9} \\
&\quad {}+ \DataV{5}{9}{15}{6}.
\end{align*}
\begin{align*}
o_{(1,1,1,1,1)}(L) &= L + 2\DataV{2}{0}{1}{0} + \DataV{3}{0}{0}{1} + \DataV{3}{0}{1}{1} \\
&\quad {}+ \DataV{3}{1}{0}{0} + \DataV{3}{1}{1}{0} + 2\DataV{3}{1}{1}{1} + \DataV{4}{0}{1}{0} \\
&\quad {}+ \DataV{4}{0}{2}{1} + \DataV{4}{1}{2}{0} + \DataV{4}{1}{3}{2} + \DataV{4}{2}{3}{1} \\
&\quad {}+ \DataV{4}{2}{4}{2} + \DataV{5}{0}{1}{0}.
\end{align*}

%% file: products/level2_square_functors.tex
\begin{align*}
\operatorname{Sym}^2\!\left(\DataV{2}{1}{0}{0}\right) &= \DataUnit + \DataV{4}{1}{0}{0} + 2\DataV{2}{1}{0}{0} + \DataV{4}{1}{1}{0} \\
&\quad {}+ \DataV{4}{1}{1}{1} + \DataV{3}{1}{1}{0} + \DataV{3}{0}{1}{0} + \DataV{4}{3}{2}{1} \\
&\quad {}+ \DataV{2}{0}{0}{1} \\
\largewedge^2\!\left(\DataV{2}{1}{0}{0}\right) &= \DataV{3}{1}{0}{0} + \DataV{3}{1}{1}{0} + \DataV{4}{2}{3}{1} + \DataV{2}{0}{1}{0} \\
&\quad {}+ L + \DataV{4}{2}{1}{0} + \DataV{3}{1}{1}{1} \\
\operatorname{Sym}^2\!\left(\DataV{2}{0}{1}{0}\right) &= \DataUnit + 2\DataV{2}{1}{0}{0} + \DataV{4}{1}{1}{0} + \DataV{4}{1}{1}{1} \\
&\quad {}+ \DataV{3}{1}{1}{0} + 2\DataV{3}{0}{1}{0} + \DataV{4}{0}{1}{0} + \DataV{4}{0}{1}{1} \\
&\quad {}+ \DataV{3}{1}{1}{1} + \DataV{4}{3}{2}{1} + \DataV{4}{2}{4}{2} + \DataV{3}{0}{1}{1} \\
&\quad {}+ 2\DataV{2}{0}{0}{1} + \DataV{4}{1}{2}{3} \\
\largewedge^2\!\left(\DataV{2}{0}{1}{0}\right) &= \DataV{3}{1}{0}{0} + \DataV{3}{1}{1}{0} + \DataV{4}{2}{3}{1} + 2\DataV{2}{0}{1}{0} \\
&\quad {}+ \DataV{4}{1}{2}{0} + \DataV{4}{1}{3}{2} + L + 2\DataV{3}{1}{1}{1} \\
&\quad {}+ \DataV{4}{2}{4}{2} + \DataV{4}{0}{2}{1} + \DataV{3}{0}{1}{1} + \DataV{3}{0}{0}{1} \\
\operatorname{Sym}^2\!\left(\DataV{2}{0}{0}{1}\right) &= \DataUnit + \DataV{2}{1}{0}{0} + \DataV{4}{1}{1}{1} + \DataV{3}{0}{1}{0} \\
&\quad {}+ \DataV{4}{0}{1}{1} + \DataV{3}{0}{1}{1} + 2\DataV{2}{0}{0}{1} + \DataV{4}{1}{2}{3} \\
&\quad {}+ \DataV{4}{0}{0}{1} \\
\largewedge^2\!\left(\DataV{2}{0}{0}{1}\right) &= \DataV{2}{0}{1}{0} + \DataV{4}{1}{3}{2} + L + \DataV{3}{1}{1}{1} \\
&\quad {}+ \DataV{3}{0}{1}{1} + \DataV{4}{0}{1}{2} + \DataV{3}{0}{0}{1} \\
\end{align*}

%% file: characters/level0.tex
\subsubsection*{Level 0 characters}
\paragraph{classical and smaller exceptional types.}
\scriptsize
\begin{longtable}{@{}>{$}l<{$} >{\raggedright\arraybackslash$}p{0.19\linewidth}<{$} >{\raggedright\arraybackslash$}p{0.19\linewidth}<{$} >{\raggedright\arraybackslash$}p{0.19\linewidth}<{$} >{\raggedright\arraybackslash$}p{0.19\linewidth}<{$}@{}}
\toprule
\text{label} & A1 & A2 & G2 & D4 \\
\midrule\endhead
\DataUnit & [0] & [0,0] & [0,0] & [0,0,0,0] \\
\bottomrule
\end{longtable}
\normalsize
\paragraph{exceptional types.}
\scriptsize
\begin{longtable}{@{}>{$}l<{$} >{\raggedright\arraybackslash$}p{0.19\linewidth}<{$} >{\raggedright\arraybackslash$}p{0.19\linewidth}<{$} >{\raggedright\arraybackslash$}p{0.19\linewidth}<{$} >{\raggedright\arraybackslash$}p{0.19\linewidth}<{$}@{}}
\toprule
\text{label} & F4 & E6 & E7 & E8 \\
\midrule\endhead
\DataUnit & [0,0,0,0] & [0,0,0,0,0,0] & [0,0,0,0,0,0,0] & [0,0,0,0,0,0,0,0] \\
\bottomrule
\end{longtable}
\normalsize

%% file: characters/level1.tex
\subsubsection*{Level 1 characters}
\paragraph{classical and smaller exceptional types.}
\scriptsize
\begin{longtable}{@{}>{$}l<{$} >{\raggedright\arraybackslash$}p{0.19\linewidth}<{$} >{\raggedright\arraybackslash$}p{0.19\linewidth}<{$} >{\raggedright\arraybackslash$}p{0.19\linewidth}<{$} >{\raggedright\arraybackslash$}p{0.19\linewidth}<{$}@{}}
\toprule
\text{label} & A1 & A2 & G2 & D4 \\
\midrule\endhead
L & [2] & [1,1] & [0,1] & [0,1,0,0] \\
\bottomrule
\end{longtable}
\normalsize
\paragraph{exceptional types.}
\scriptsize
\begin{longtable}{@{}>{$}l<{$} >{\raggedright\arraybackslash$}p{0.19\linewidth}<{$} >{\raggedright\arraybackslash$}p{0.19\linewidth}<{$} >{\raggedright\arraybackslash$}p{0.19\linewidth}<{$} >{\raggedright\arraybackslash$}p{0.19\linewidth}<{$}@{}}
\toprule
\text{label} & F4 & E6 & E7 & E8 \\
\midrule\endhead
L & [1,0,0,0] & [0,1,0,0,0,0] & [1,0,0,0,0,0,0] & [0,0,0,0,0,0,0,1] \\
\bottomrule
\end{longtable}
\normalsize

%% file: characters/level2.tex
\subsubsection*{Level 2 characters}
\paragraph{classical and smaller exceptional types.}
\scriptsize
\begin{longtable}{@{}>{$}l<{$} >{\raggedright\arraybackslash$}p{0.19\linewidth}<{$} >{\raggedright\arraybackslash$}p{0.19\linewidth}<{$} >{\raggedright\arraybackslash$}p{0.19\linewidth}<{$} >{\raggedright\arraybackslash$}p{0.19\linewidth}<{$}@{}}
\toprule
\text{label} & A1 & A2 & G2 & D4 \\
\midrule\endhead
\DataV{2}{1}{0}{0} & 0 & [1,1] & [2,0] & [0,0,0,2]\allowbreak +\nobreak\,[0,0,2,0]\allowbreak +\nobreak\,[2,0,0,0] \\
\DataV{2}{0}{1}{0} & 0 & [0,3]\allowbreak +\nobreak\,[3,0] & [3,0] & [1,0,1,1] \\
\DataV{2}{0}{0}{1} & [4] & [2,2] & [0,2] & [0,2,0,0] \\
\bottomrule
\end{longtable}
\normalsize
\paragraph{exceptional types.}
\scriptsize
\begin{longtable}{@{}>{$}l<{$} >{\raggedright\arraybackslash$}p{0.19\linewidth}<{$} >{\raggedright\arraybackslash$}p{0.19\linewidth}<{$} >{\raggedright\arraybackslash$}p{0.19\linewidth}<{$} >{\raggedright\arraybackslash$}p{0.19\linewidth}<{$}@{}}
\toprule
\text{label} & F4 & E6 & E7 & E8 \\
\midrule\endhead
\DataV{2}{1}{0}{0} & [0,0,0,2] & [1,0,0,0,0,1] & [0,0,0,0,0,1,0] & [1,0,0,0,0,0,0,0] \\
\DataV{2}{0}{1}{0} & [0,1,0,0] & [0,0,0,1,0,0] & [0,0,1,0,0,0,0] & [0,0,0,0,0,0,1,0] \\
\DataV{2}{0}{0}{1} & [2,0,0,0] & [0,2,0,0,0,0] & [2,0,0,0,0,0,0] & [0,0,0,0,0,0,0,2] \\
\bottomrule
\end{longtable}
\normalsize

%% file: characters/level3.tex
\subsubsection*{Level 3 characters}
\paragraph{classical and smaller exceptional types.}
\scriptsize
\begin{longtable}{@{}>{$}l<{$} >{\raggedright\arraybackslash$}p{0.19\linewidth}<{$} >{\raggedright\arraybackslash$}p{0.19\linewidth}<{$} >{\raggedright\arraybackslash$}p{0.19\linewidth}<{$} >{\raggedright\arraybackslash$}p{0.19\linewidth}<{$}@{}}
\toprule
\text{label} & A1 & A2 & G2 & D4 \\
\midrule\endhead
\DataV{3}{1}{0}{0} & 0 & [0,0] & [1,0] & 2[0,1,0,0] \\
\DataV{3}{0}{1}{0} & -[4] & 0 & [4,0] & [0,0,2,2]\allowbreak +\nobreak\,[2,0,0,2]\allowbreak +\nobreak\,[2,0,2,0] \\
\DataV{3}{0}{0}{1} & [6] & [3,3] & [0,3] & [0,3,0,0] \\
\DataV{3}{0}{1}{1} & 0 & [1,4]\allowbreak +\nobreak\,[4,1] & [3,1] & [1,1,1,1] \\
\DataV{3}{1}{1}{0} & -[2] & 0 & [1,1] & 2[1,0,1,1] \\
\DataV{3}{1}{1}{1} & 0 & [2,2] & [2,1] & [0,1,0,2]\allowbreak +\nobreak\,[0,1,2,0]\allowbreak +\nobreak\,[2,1,0,0] \\
\bottomrule
\end{longtable}
\normalsize
\paragraph{exceptional types.}
\scriptsize
\begin{longtable}{@{}>{$}l<{$} >{\raggedright\arraybackslash$}p{0.19\linewidth}<{$} >{\raggedright\arraybackslash$}p{0.19\linewidth}<{$} >{\raggedright\arraybackslash$}p{0.19\linewidth}<{$} >{\raggedright\arraybackslash$}p{0.19\linewidth}<{$}@{}}
\toprule
\text{label} & F4 & E6 & E7 & E8 \\
\midrule\endhead
\DataV{3}{1}{0}{0} & [0,0,1,0] & [1,0,0,0,0,1] & [0,0,0,0,0,0,2] & 0 \\
\DataV{3}{0}{1}{0} & [0,0,2,0] & [0,0,1,0,1,0] & [0,0,0,1,0,0,0] & [0,0,0,0,0,1,0,0] \\
\DataV{3}{0}{0}{1} & [3,0,0,0] & [0,3,0,0,0,0] & [3,0,0,0,0,0,0] & [0,0,0,0,0,0,0,3] \\
\DataV{3}{0}{1}{1} & [1,1,0,0] & [0,1,0,1,0,0] & [1,0,1,0,0,0,0] & [0,0,0,0,0,0,1,1] \\
\DataV{3}{1}{1}{0} & [0,0,1,1] & [0,0,0,0,1,1]\allowbreak +\nobreak\,[1,0,1,0,0,0] & [0,1,0,0,0,0,1] & [0,1,0,0,0,0,0,0] \\
\DataV{3}{1}{1}{1} & [1,0,0,2] & [1,1,0,0,0,1] & [1,0,0,0,0,1,0] & [1,0,0,0,0,0,0,1] \\
\bottomrule
\end{longtable}
\normalsize

%% file: characters/level4.tex
\subsubsection*{Level 4 characters}
\paragraph{classical and smaller exceptional types.}
\scriptsize
\begin{longtable}{@{}>{$}l<{$} >{\raggedright\arraybackslash$}p{0.19\linewidth}<{$} >{\raggedright\arraybackslash$}p{0.19\linewidth}<{$} >{\raggedright\arraybackslash$}p{0.19\linewidth}<{$} >{\raggedright\arraybackslash$}p{0.19\linewidth}<{$}@{}}
\toprule
\text{label} & A1 & A2 & G2 & D4 \\
\midrule\endhead
\DataV{4}{1}{0}{0} & 0 & 0 & 0 & 2[0,0,0,0] \\
\DataV{4}{0}{0}{1} & [8] & [4,4] & [0,4] & [0,4,0,0] \\
\DataV{4}{0}{1}{0} & 0 & -[1,4]\allowbreak -\nobreak\,[4,1] & 0 & [1,0,1,3]\allowbreak +\nobreak\,[1,0,3,1]\allowbreak +\nobreak\,[3,0,1,1] \\
\DataV{4}{0}{1}{1} & 0 & [0,6]\allowbreak +\nobreak\,[6,0] & [6,0] & [2,0,2,2] \\
\DataV{4}{1}{1}{0} & -[0] & -[1,1] & 0 & [1,0,1,1] \\
\DataV{4}{1}{1}{1} & 0 & 0 & 0 & [0,0,0,4]\allowbreak +\nobreak\,[0,0,4,0]\allowbreak +\nobreak\,[4,0,0,0] \\
\DataV{4}{1}{2}{0} & 0 & -[0,3]\allowbreak -\nobreak\,[3,0] & 0 & [0,1,0,2]\allowbreak +\nobreak\,[0,1,2,0]\allowbreak +\nobreak\,[2,1,0,0] \\
\DataV{4}{2}{1}{0} & 0 & -[0,0] & 0 & [0,0,0,2]\allowbreak +\nobreak\,[0,0,2,0]\allowbreak +\nobreak\,[2,0,0,0] \\
\DataV{4}{0}{2}{1} & -[6] & 0 & [4,1] & [0,1,2,2]\allowbreak +\nobreak\,[2,1,0,2]\allowbreak +\nobreak\,[2,1,2,0] \\
\DataV{4}{0}{1}{2} & 0 & [2,5]\allowbreak +\nobreak\,[5,2] & [3,2] & [1,2,1,1] \\
\DataV{4}{1}{3}{2} & 0 & 0 & [5,0] & [1,0,1,3]\allowbreak +\nobreak\,[1,0,3,1]\allowbreak +\nobreak\,[3,0,1,1] \\
\DataV{4}{1}{2}{3} & 0 & [3,3] & [2,2] & [0,2,0,2]\allowbreak +\nobreak\,[0,2,2,0]\allowbreak +\nobreak\,[2,2,0,0] \\
\DataV{4}{3}{2}{1} & [2] & 0 & 0 & 2[0,2,0,0] \\
\DataV{4}{2}{3}{1} & 0 & -[2,2] & 0 & [0,0,2,2]\allowbreak +\nobreak\,[2,0,0,2]\allowbreak +\nobreak\,[2,0,2,0] \\
\DataV{4}{2}{4}{2} & 0 & 0 & [1,2] & 2[1,1,1,1] \\
\bottomrule
\end{longtable}
\normalsize
\paragraph{exceptional types.}
\scriptsize
\begin{longtable}{@{}>{$}l<{$} >{\raggedright\arraybackslash$}p{0.19\linewidth}<{$} >{\raggedright\arraybackslash$}p{0.19\linewidth}<{$} >{\raggedright\arraybackslash$}p{0.19\linewidth}<{$} >{\raggedright\arraybackslash$}p{0.19\linewidth}<{$}@{}}
\toprule
\text{label} & F4 & E6 & E7 & E8 \\
\midrule\endhead
\DataV{4}{1}{0}{0} & [0,0,0,1] & [0,1,0,0,0,0] & 0 & -[1,0,0,0,0,0,0,0] \\
\DataV{4}{0}{0}{1} & [4,0,0,0] & [0,4,0,0,0,0] & [4,0,0,0,0,0,0] & [0,0,0,0,0,0,0,4] \\
\DataV{4}{0}{1}{0} & [0,0,2,1] & [0,0,2,0,0,1]\allowbreak +\nobreak\,[1,0,0,0,2,0] & [0,1,0,0,1,0,0] & [0,0,0,0,1,0,0,0] \\
\DataV{4}{0}{1}{1} & [0,2,0,0] & [0,0,0,2,0,0] & [0,0,2,0,0,0,0] & [0,0,0,0,0,0,2,0] \\
\DataV{4}{1}{1}{0} & [0,0,0,3] & [0,0,0,0,0,3]\allowbreak +\nobreak\,[3,0,0,0,0,0] & 0 & 0 \\
\DataV{4}{1}{1}{1} & [0,0,0,4] & [2,0,0,0,0,2] & [0,0,0,0,0,2,0] & [2,0,0,0,0,0,0,0] \\
\DataV{4}{1}{2}{0} & [0,1,0,1] & [0,0,1,0,1,0] & [0,2,0,0,0,0,0] & 0 \\
\DataV{4}{2}{1}{0} & [1,0,0,1] & [0,0,0,1,0,0] & 0 & -[0,1,0,0,0,0,0,0] \\
\DataV{4}{0}{2}{1} & [1,0,2,0] & [0,1,1,0,1,0] & [1,0,0,1,0,0,0] & [0,0,0,0,0,1,0,1] \\
\DataV{4}{0}{1}{2} & [2,1,0,0] & [0,2,0,1,0,0] & [2,0,1,0,0,0,0] & [0,0,0,0,0,0,1,2] \\
\DataV{4}{1}{3}{2} & [0,1,0,2] & [1,0,0,1,0,1] & [0,0,1,0,0,1,0] & [1,0,0,0,0,0,1,0] \\
\DataV{4}{1}{2}{3} & [2,0,0,2] & [1,2,0,0,0,1] & [2,0,0,0,0,1,0] & [1,0,0,0,0,0,0,2] \\
\DataV{4}{3}{2}{1} & [1,0,1,0] & [1,1,0,0,0,1] & [1,0,0,0,0,0,2] & 0 \\
\DataV{4}{2}{3}{1} & [0,0,1,2] & [0,0,1,0,0,2]\allowbreak +\nobreak\,[2,0,0,0,1,0] & [0,0,0,0,1,0,1] & [0,0,1,0,0,0,0,0] \\
\DataV{4}{2}{4}{2} & [1,0,1,1] & [0,1,0,0,1,1]\allowbreak +\nobreak\,[1,1,1,0,0,0] & [1,1,0,0,0,0,1] & [0,1,0,0,0,0,0,1] \\
\bottomrule
\end{longtable}
\normalsize

%% file: characters/level5.tex
\subsubsection*{Level 5 characters}
\paragraph{classical and smaller exceptional types.}
\scriptsize
\begin{longtable}{@{}>{$}l<{$} >{\raggedright\arraybackslash$}p{0.19\linewidth}<{$} >{\raggedright\arraybackslash$}p{0.19\linewidth}<{$} >{\raggedright\arraybackslash$}p{0.19\linewidth}<{$} >{\raggedright\arraybackslash$}p{0.19\linewidth}<{$}@{}}
\toprule
\text{label} & A1 & A2 & G2 & D4 \\
\midrule\endhead
\DataV{5}{1}{0}{0} & 0 & 0 & 0 & 0 \\
\DataV{5}{0}{0}{1} & [10] & [5,5] & [0,5] & [0,5,0,0] \\
\DataV{5}{0}{1}{0} & 0 & -[3,3] & -[4,1] & [0,1,0,4]\allowbreak +\nobreak\,[0,1,4,0]\allowbreak +\nobreak\,[2,0,2,2]\allowbreak +\nobreak\,[4,1,0,0] \\
\DataV{5}{0}{1}{3} & 0 & [3,6]\allowbreak +\nobreak\,[6,3] & [3,3] & [1,3,1,1] \\
\DataV{5}{1}{3}{0} & 0 & 0 & 0 & 2[0,0,0,4]\allowbreak +\nobreak\,2[0,0,4,0]\allowbreak +\nobreak\,2[4,0,0,0] \\
\DataV{5}{0}{3}{1} & 0 & -[2,5]\allowbreak -\nobreak\,[5,2] & 0 & [1,1,1,3]\allowbreak +\nobreak\,[1,1,3,1]\allowbreak +\nobreak\,[3,1,1,1] \\
\DataV{5}{3}{1}{0} & 0 & 0 & 0 & -2[0,0,0,0] \\
\DataV{5}{0}{2}{3} & 0 & [1,7]\allowbreak +\nobreak\,[7,1] & [6,1] & [2,1,2,2] \\
\DataV{5}{0}{3}{2} & 0 & 0 & [7,0] & [1,0,3,3]\allowbreak +\nobreak\,[3,0,1,3]\allowbreak +\nobreak\,[3,0,3,1] \\
\DataV{5}{2}{3}{0} & 0 & -[1,1] & 0 & 3[0,2,0,0] \\
\DataV{5}{3}{2}{0} & 0 & 0 & 0 & [0,1,0,0] \\
\DataV{5}{0}{3}{3} & -[8] & 0 & [4,2] & [0,2,2,2]\allowbreak +\nobreak\,[2,2,0,2]\allowbreak +\nobreak\,[2,2,2,0] \\
\DataV{5}{3}{3}{0} & 0 & [0,0] & 0 & 0 \\
\DataV{5}{1}{6}{3} & 0 & -[0,6]\allowbreak -\nobreak\,[6,0] & 0 & [0,0,2,4]\allowbreak +\nobreak\,[0,0,4,2]\allowbreak +\nobreak\,[2,0,0,4]\allowbreak +\nobreak\,[2,0,4,0]\allowbreak +\nobreak\,[4,0,0,2]\allowbreak +\nobreak\,[4,0,2,0] \\
\DataV{5}{3}{6}{1} & 0 & -[2,2] & -[4,0] & -[0,0,0,4]\allowbreak -\nobreak\,[0,0,4,0]\allowbreak -\nobreak\,[4,0,0,0] \\
\DataV{5}{1}{3}{6} & 0 & [4,4] & [2,3] & [0,3,0,2]\allowbreak +\nobreak\,[0,3,2,0]\allowbreak +\nobreak\,[2,3,0,0] \\
\DataV{5}{6}{3}{1} & 0 & 0 & -[1,0] & -[0,0,0,2]\allowbreak -\nobreak\,[0,0,2,0]\allowbreak -\nobreak\,[2,0,0,0] \\
\DataV{5}{3}{6}{6} & 0 & 0 & 0 & [0,1,0,4]\allowbreak +\nobreak\,[0,1,4,0]\allowbreak +\nobreak\,[4,1,0,0] \\
\DataV{5}{6}{6}{3} & 0 & 0 & -[2,1] & -[0,0,2,2]\allowbreak -\nobreak\,[2,0,0,2]\allowbreak -\nobreak\,[2,0,2,0] \\
\DataV{5}{3}{9}{3} & 0 & 0 & 0 & [0,2,0,2]\allowbreak +\nobreak\,[0,2,2,0]\allowbreak +\nobreak\,[2,2,0,0] \\
\DataV{5}{2}{9}{9} & 0 & 0 & [5,1] & [1,1,1,3]\allowbreak +\nobreak\,[1,1,3,1]\allowbreak +\nobreak\,[3,1,1,1] \\
\DataV{5}{6}{8}{6} & 0 & 0 & 0 & 2[0,3,0,0] \\
\DataV{5}{3}{9}{8} & 0 & 0 & [1,3] & 2[1,2,1,1] \\
\DataV{5}{8}{9}{3} & 0 & [0,3]\allowbreak +\nobreak\,[3,0] & 0 & 0 \\
\DataV{5}{3}{11}{6} & [6] & 0 & 0 & 2[2,0,2,2] \\
\DataV{5}{9}{9}{2} & [0] & 0 & -[1,1] & -2[0,2,0,0] \\
\DataV{5}{6}{11}{3} & 0 & 0 & 0 & [1,1,1,1] \\
\DataV{5}{6}{12}{6} & 0 & 0 & -[5,0] & 0 \\
\DataV{5}{6}{15}{9} & 0 & -[3,3] & 0 & [0,1,2,2]\allowbreak +\nobreak\,[2,1,0,2]\allowbreak +\nobreak\,[2,1,2,0] \\
\DataV{5}{9}{15}{6} & [4] & 0 & -[1,2] & 0 \\
\bottomrule
\end{longtable}
\normalsize
\paragraph{exceptional types.}
\scriptsize
\begin{longtable}{@{}>{$}l<{$} >{\raggedright\arraybackslash$}p{0.19\linewidth}<{$} >{\raggedright\arraybackslash$}p{0.19\linewidth}<{$} >{\raggedright\arraybackslash$}p{0.19\linewidth}<{$} >{\raggedright\arraybackslash$}p{0.19\linewidth}<{$}@{}}
\toprule
\text{label} & F4 & E6 & E7 & E8 \\
\midrule\endhead
\DataV{5}{1}{0}{0} & 0 & [0,0,0,0,0,0] & 0 & -[0,0,0,0,0,0,0,1] \\
\DataV{5}{0}{0}{1} & [5,0,0,0] & [0,5,0,0,0,0] & [5,0,0,0,0,0,0] & [0,0,0,0,0,0,0,5] \\
\DataV{5}{0}{1}{0} & [0,0,3,0]\allowbreak +\nobreak\,[0,1,0,3] & [0,0,0,0,3,0]\allowbreak +\nobreak\,[0,0,3,0,0,0]\allowbreak +\nobreak\,[1,0,1,0,1,1] & [0,0,0,0,2,0,0]\allowbreak +\nobreak\,[0,2,0,0,0,1,0] & [0,0,0,1,0,0,0,0] \\
\DataV{5}{0}{1}{3} & [3,1,0,0] & [0,3,0,1,0,0] & [3,0,1,0,0,0,0] & [0,0,0,0,0,0,1,3] \\
\DataV{5}{1}{3}{0} & 0 & -[0,1,0,0,0,3]\allowbreak -\nobreak\,[3,1,0,0,0,0] & -[0,0,1,0,0,0,2] & -[0,0,0,0,1,0,0,0] \\
\DataV{5}{0}{3}{1} & [1,0,2,1] & [0,1,2,0,0,1]\allowbreak +\nobreak\,[1,1,0,0,2,0] & [1,1,0,0,1,0,0] & [0,0,0,0,1,0,0,1] \\
\DataV{5}{3}{1}{0} & 0 & 0 & -[0,0,0,0,0,0,2] & -[0,0,0,0,0,0,1,0] \\
\DataV{5}{0}{2}{3} & [1,2,0,0] & [0,1,0,2,0,0] & [1,0,2,0,0,0,0] & [0,0,0,0,0,0,2,1] \\
\DataV{5}{0}{3}{2} & [0,1,2,0] & [0,0,1,1,1,0] & [0,0,1,1,0,0,0] & [0,0,0,0,0,1,1,0] \\
\DataV{5}{2}{3}{0} & 0 & -[0,0,1,0,1,0] & -[0,0,0,0,1,0,1] & 0 \\
\DataV{5}{3}{2}{0} & 0 & -[0,0,0,1,0,0] & -[0,1,0,0,0,0,1] & 0 \\
\DataV{5}{0}{3}{3} & [2,0,2,0] & [0,2,1,0,1,0] & [2,0,0,1,0,0,0] & [0,0,0,0,0,1,0,2] \\
\DataV{5}{3}{3}{0} & 0 & -[0,0,0,0,0,3]\allowbreak -\nobreak\,[3,0,0,0,0,0] & -[1,0,0,0,0,0,2] & -[0,0,0,0,0,1,0,0] \\
\DataV{5}{1}{6}{3} & [0,0,2,2] & [1,0,1,0,1,1] & [0,0,0,1,0,1,0] & [1,0,0,0,0,1,0,0] \\
\DataV{5}{3}{6}{1} & 0 & 0 & 0 & 0 \\
\DataV{5}{1}{3}{6} & [3,0,0,2] & [1,3,0,0,0,1] & [3,0,0,0,0,1,0] & [1,0,0,0,0,0,0,3] \\
\DataV{5}{6}{3}{1} & 0 & [0,2,0,0,0,0] & 0 & -[1,0,0,0,0,0,0,1] \\
\DataV{5}{3}{6}{6} & [1,0,0,4] & [2,1,0,0,0,2] & [1,0,0,0,0,2,0] & [2,0,0,0,0,0,0,1] \\
\DataV{5}{6}{6}{3} & 0 & [2,0,0,0,0,2] & [0,0,0,0,0,1,2] & 0 \\
\DataV{5}{3}{9}{3} & [1,1,0,1] & [0,1,1,0,1,0] & [1,2,0,0,0,0,0] & 0 \\
\DataV{5}{2}{9}{9} & [1,1,0,2] & [1,1,0,1,0,1] & [1,0,1,0,0,1,0] & [1,0,0,0,0,0,1,1] \\
\DataV{5}{6}{8}{6} & [2,0,1,0] & [1,2,0,0,0,1] & [2,0,0,0,0,0,2] & 0 \\
\DataV{5}{3}{9}{8} & [2,0,1,1] & [0,2,0,0,1,1]\allowbreak +\nobreak\,[1,2,1,0,0,0] & [2,1,0,0,0,0,1] & [0,1,0,0,0,0,0,2] \\
\DataV{5}{8}{9}{3} & [2,0,0,1] & [0,1,0,1,0,0] & 0 & -[0,1,0,0,0,0,0,1] \\
\DataV{5}{3}{11}{6} & [0,1,1,1] & [0,0,0,1,1,1]\allowbreak +\nobreak\,[1,0,1,1,0,0] & [0,1,1,0,0,0,1] & [0,1,0,0,0,0,1,0] \\
\DataV{5}{9}{9}{2} & 0 & 0 & -[0,2,0,0,0,0,0] & -[0,0,1,0,0,0,0,0] \\
\DataV{5}{6}{11}{3} & [1,0,0,3] & [0,1,0,0,0,3]\allowbreak +\nobreak\,[3,1,0,0,0,0] & 0 & 0 \\
\DataV{5}{6}{12}{6} & [0,0,1,3] & [1,0,0,0,1,2]\allowbreak +\nobreak\,[2,0,1,0,0,1] & [0,1,0,0,0,1,1] & [1,1,0,0,0,0,0,0] \\
\DataV{5}{6}{15}{9} & [1,0,1,2] & [0,1,1,0,0,2]\allowbreak +\nobreak\,[2,1,0,0,1,0] & [1,0,0,0,1,0,1] & [0,0,1,0,0,0,0,1] \\
\DataV{5}{9}{15}{6} & [0,1,1,0] & [1,0,0,1,0,1] & [0,0,1,0,0,0,2] & 0 \\
\bottomrule
\end{longtable}
\normalsize